\documentclass{amsart}
  \newtheorem{thm}{Theorem}[section]
  \newtheorem{lem}[thm]{Lemma}
  \newtheorem{defn}[thm]{Definition}
  \usepackage{upgreek} 
  \usepackage[linesnumbered]{algorithm2e}
  \usepackage{graphicx} 

 \title[Simulation and Control of a Nonsmooth CHNS System]{Simulation and Control of a Nonsmooth Cahn-Hilliard Navier-Stokes System with Variable Fluid Densities}

 \author[Gr\"a\ss{}le]{Carmen Gr\"a\ss{}le}
 
   \address{%
   University of Hamburg\\
   Bundesstr.\ 55\\
   20146 Hamburg\\
   Germany}
  
  \email{carmen.graessle@uni-hamburg.de}

 \author[Hinterm\"uller]{Michael Hinterm\"uller}
 \address{Weierstra\ss{}-Institut\\
 Mohrenstr.\ 39\\
 10117 Berlin\\
 Germany}
 \email{michael.hintermueller@wias-berlin.de}
 \author[Hinze]{Michael Hinze}
 \address{University of Koblenz-Landau\\
 Universit\"atsstr.\ 1\\
 56070 Koblenz\\
 Germany}
 \email{hinze@uni-koblenz.de}
 \author[Keil]{Tobias Keil}
 \address{Weierstra\ss{}-Institut\\
 Mohrenstr.\ 39\\
 10117 Berlin
 Germany}
 \email{tobias.keil@wias-berlin.de}

\date{}
 
 \subjclass{35B65, 35J87, 35K55, 65K15, 49K20}
 
 \keywords{Two phase flow, POD model order reduction, adaptivity, nonsmooth systems, mathematical programming with equilibrium constraints, optimal control}

\begin{document}
\maketitle

 \begin{abstract}
 We are concerned with the simulation and control of a two phase flow model governed by a coupled Cahn-Hilliard Navier-Stokes system involving a nonsmooth energy potential.
 We establish the existence of optimal solutions and present two distinct approaches to derive suitable stationarity conditions for the bilevel problem, namely C- and strong stationarity.
 Moreover, we demonstrate the numerical realization of these concepts at the hands of two adaptive solution algorithms relying on a specifically developed goal-oriented error estimator.
 In addition, we present a model order reduction approach using proper orthogonal decomposition (POD-MOR) in order to replace high-fidelity models by low order surrogates. In particular, we combine POD with space-adapted snapshots and address the challenges which are the consideration of snapshots with different spatial resolutions and the conservation of a solenoidal property.
 \end{abstract}
%

 \section{Introduction}
 
 We consider the simulation and control for multiphase flows governed by a Cahn-Hilliard Navier-Stokes system with nonsmooth homogeneous free energy densities utilizing a diffuse interface approach. The free energy is a double-obstacle potential according to \cite{BE91}. The resulting problem belongs to the class of mathematical programs with equilibrium constraints (MPECs) in function space.
 
 Even in finite dimensions, this problem class is well-known for its constraint degeneracy \cite{Luo1996,Outrata1998}. Due to the presence of the variational inequality constraint, classical constraint qualifications (see, e.g., \cite{Zowe1979}) fail which prevents the application of Karush-Kuhn-Tucker (KKT) theory in Banach space for the first-order characterization of an optimal solution by (Lagrange) multipliers. As a result, stationarity conditions for this problem class are no longer unique (in contrast to KKT conditions); compare \cite{Hintermuller2009,Hintermuller2014} in function space and, e.g., \cite{Scheel2000} in finite dimensions. They rather depend on the underlying problem structure and/or on the chosen analytical approach.

 The simulation of two-phase flows with matched densities is rather well understood in the literature, see e.g. \cite{Hohenberg1977}.
 In contrast, there exist different approaches to model the case of fluids with non-matched densities.
 These range from quasi-incompressible models with non-divergence free velocity fields, see e.g. \cite{Lowengrub1998}, to possibly thermodynamically inconsistent models with solenoidal fluid velocities, cf. \cite{Ding2007}.
 In this work, we study the incompressible and thermodynamically consistent model presented in \cite{AGG}.
 We refer to \cite{n15_2012_AlandVoigt_bubble_benchmark,Boyer2002,Boyer2004,Gal2010,GHK} for additional analytical and numerical results for some of these models.

 Stable numerical schemes for the thermodynamically consistent diffuse interface model according to \cite{AGG} are developed in \cite{GHK,GK}. A fully integrated adaptive finite element approach for the numerical treatment of the Cahn-Hilliard system with a nonsmooth free energy is developed in \cite{HHT}. This approach is extended in \cite{HHK12} to a fully practical adaptive solver for the coupled Cahn-Hilliard Navier-Stokes system.
 
 While there are numerous publications concerning the optimal control of the phase separation process itself, i.e. the distinct Cahn--Hilliard system,
 see e.g. \cite{BE91,Colli2014,Elliott1986,HHT,Hintermueller2012,Yong1991}, there has been considerably less research on the control of the Cahn--Hilliard--Navier--Stokes system.
 Some of the few publications in this field address the case of matched densities and a non-smooth homogeneous free energy density (double-obstacle potential), see \cite{Hintermuller2014a,Hintermueller2014}.
 We also mention the recent  {articles} \cite{Frigeri2014a}, which treats the control of a nonlocal Cahn--Hilliard--Navier--Stokes 
 system in two dimensions, \cite{Tachim2015}, and \cite{GHK19}, which includes numerical convergence results for the optimal control of the model developed in \cite{GHK}.  

 From a numerical point of view, the simulation and especially the optimal control of the coupled Cahn-Hilliard Navier-Stokes system are challenging tasks in regards to the computational times and the storage effort. For this reason, we apply model order reduction using Proper Orthogonal Decomposition (POD-MOR) in order to replace the high-fidelity models by low-order surrogates. We follow a simulation-based approach according to \cite{Sir}, where the snapshots are generated by finite element simulations of the system. In particular, we utilize space-adapted snapshots which leads to the challenge that in a discrete formulation the snapshots are vectors of different lengths due to the different spatial resolutions. A consideration of the problem setting from an infinite-dimensional view according to \cite{GH18} allows the combination of POD with spatially adapted snapshots. Moreover, we utilize a Moreau-Yosida regularization in the Cahn-Hilliard system and observe that the accuracy of the reduced-order model depends on the smoothness of the approximated object. Finally, we consider POD-MOR for the Navier-Stokes part. The use of space-adapted finite elements has the consequence that a weak divergence-free property only holds in the current adapted finite element space. In order to guarantee stability of the resulting reduced-order model, in \cite{GHLU18} two solution approaches are proposed.

 Regarding physical applications, we point out that the CHNS system is used to model a variety of situations.
 These range from the aforementioned solidification process of liquid
 metal alloys, cf. \cite{Eckert2013}, or the simulation of bubble dynamics, as in Taylor flows \cite{Aland2013}, or pinch-offs of liquid-liquid jets \cite{Kim2004},
 to the formation of polymeric membranes \cite{Zhou2006} or proteins
 crystallization, see e.g. \cite{Kim2004a} and references within.
 Furthermore, the model can be easily adapted to include the effects of surfactants such as colloid
 particles at fluid-fluid interfaces in gels and emulsions used in food, pharmaceutical, cosmetic, or petroleum industries \cite{Aland2012,Praetorius2013}.

 The paper is organized as follows. After introducing the problem setting in Section \ref{p13sec:problemsetting}, we formulate the associated optimal control problem with respect to a semi-discrete system in Section \ref{p13sec:OC}.
 We proceed by securing the existence of global solutions and characterizing these solutions via suitable stationarity conditions in the Sections \ref{p13sec:exist}, \ref{p13sec:Cstat} and \ref{p13sec:Sstat}.
 We present two distinct numerical solution algorithms based on our analytical results in Section \ref{p13sec:penalgo} and \ref{p13sec:bundle} and incorporate an adaptive mesh refinement technique relying on a goal-oriented error estimator of Section \ref{p13sec:adapt}. 
 In Section \ref{p13sec:mor}, we focus on model order reduction with Proper Orthogonal Decomposition. The POD method in Hilbert spaces is explained in Section \ref{p13sec:PODHilb} and comprises the case of space-adapted snapshots. In Section \ref{p13sec:PODCH} we derive a POD reduced-order model for the Cahn-Hilliard system and provide a numerical example in Section \ref{p13sec:numericspodch}. Moreover, in Section \ref{p13sec:PODNaSt} we consider POD-MOR with space-adapted snapshots for the Navier-Stokes equations.
 We conclude this article with a brief outlook on associated future research topics in Section \ref{p13sec:out}.
 
 \section{Problem setting}\label{p13sec:problemsetting}
 
 Let us specify the problem setting. We denote by $\Omega$ an open bounded domain with Lipschitz boundary $\partial \Omega$ and $T>0$ is a given end time. We are concerned with the coupled Cahn-Hilliard Navier-Stokes (CHNS) system according to \cite{AGG} given by
  \begin{subequations}\label{p13eq:CHNS}
  \begin{alignat}{4}
 \partial_t(\rho(\varphi) v)+\textnormal{div}( v\otimes\rho(\varphi) v)-\textnormal{div}(2\eta(\varphi)\epsilon(v))+\nabla p & & & \nonumber \\
  +\textnormal{div}(v\otimes -\frac{\widehat{\rho}_2-\widehat{\rho}_1}{2}m(\varphi)\nabla\mu) -\mu\nabla\varphi & = && \; 0 &&  \quad \text{in } (0,T)\times \Omega, \label{p13eq:NaSt1}\\
  \textnormal{div}v & = && \; 0 && \quad\text{in } (0,T)\times \Omega, \label{p13eq:NaSt2}\\
  \partial_t\varphi +v\nabla\varphi  -\textnormal{div}(m(\varphi)\nabla\mu) & = && \; 0 && \quad\text{in } (0,T)\times \Omega, \label{p13eq:CH1} \\
  -\sigma \epsilon\Delta\varphi + \frac{\sigma}{\epsilon}(\partial \Psi_0(\varphi) -\kappa\varphi ) -\mu &\ni && \; 0 && \quad\text{in } (0,T)\times \Omega, \label{p13eq:CH2} \\
   v=\partial_n\varphi
  =\partial_n\mu & = && \; 0  && \quad\text{on } (0,T)\times \partial\Omega,  \label{p13eq:bdryCH}\\
  v(0,\cdot) & = && \; v_a && \quad\text{in } \Omega, \label{p13eq:initialNaSt}\\
  \varphi(0,\cdot) & = && \; \varphi_a && \quad\text{in } \Omega. \label{p13eq:initialCH}
 \end{alignat}\label{p13CHNS}
 \end{subequations}
 
 We denote by $v$ the velocity and by $p$ the pressure of the fluid which is governed by the Navier-Stokes equations \eqref{p13eq:NaSt1}-\eqref{p13eq:NaSt2}. The density $\rho$ depends on the order parameter $\varphi$ given by the Cahn-Hilliard equations \eqref{p13eq:CH1}-\eqref{p13eq:CH2} via
 \begin{align}
  \rho(\varphi)=\frac{\rho_1+\rho_2}{2}+\frac{\rho_2-\rho_1}{2}\varphi.\label{p13denscon} 
 \end{align}
 The mobility $m$ and the viscosity $\eta$ are variable and depend on the phase field $\varphi$.
 By $\mu$ we denote the chemical potential. The surface tension $\sigma>0$, the interface parameter $\epsilon>0$ and the parameter $\kappa>0$ are given constants.
 Furthermore, initial conditions $v_a$ and $\varphi_a$ for the velocity and phase field are given, respectively.
 By $\Psi_0$ we denote the convex part of the free energy potential $\Psi(\varphi):=(\Psi_0(\varphi)-\frac{\kappa}{2} \varphi^2)$.
 and accounts for the restriction of the phase field variable to stay in the physically meaningful range of $[-1,1]$.
 Depending on the underlying applications, there exist different modeling choices for $\Psi_0$.
 In this article, we focus on the double-obstacle potential introduced in \eqref{p13dob}.
 Possible other choices include the double-well potential $\Psi(\varphi)=\frac{\kappa}{2}(1-\varphi^2)^2$ and the logarithmic potential $\Psi(\varphi)=(1+\varphi)\ln (1+\varphi)+(1-\varphi)\ln (1-\varphi) -\frac{\kappa}{2}\varphi^2 $.

 An important property of the above CHNS system is its thermodynamical consistency. 
It is possible to derive a (dissipative) energy estimate by testing \eqref{p13eq:NaSt1},\eqref{p13eq:NaSt2},\eqref{p13eq:CH1}, and \eqref{p13eq:CH2} with  $v$, $p$, $\mu$, and $\partial_t\varphi$, which yields
\begin{align}
 \partial_t E(v,\varphi) + 2 \int_\Omega \eta(\varphi) |\epsilon(v)|^2 dx +  \int_\Omega m(\varphi) |\nabla \mu|^2 dx \leq 0,\label{p13energyIE}
\end{align}
where the total energy $E$ is given by the sum of the kinetic and the potential energy, i.e.
\begin{align}
 E(v,\varphi)=\int_\Omega \rho(\varphi)\frac{|v|^2}{2} dx + \frac{\sigma\epsilon}{2}\int_\Omega \frac{|\nabla\varphi|^2}{2} dx + \frac{\sigma}{\epsilon}\Psi(\varphi).\label{p13totenc}
\end{align}
Besides mirroring the physical property that the total energy of a closed system is non-increasing, inequality \eqref{p13energyIE} also serves as a very valuable analytical tool, e.g., to secure the boundedness of solutions to \eqref{p13CHNS}.

\section{Optimal control of the semi-discrete CHNS system}
\label{p13Chap:OC}

In the following, we study the optimal control of a semi-discrete variant of the Cahn-Hilliard Navier-Stokes system \eqref{p13CHNS}, where the free energy density is related to the double-obstacle potential, see \eqref{p13dob} below.
This yields an optimal control problem for a family of coupled systems in each time {instant} of a variational inequality of fourth order and the Navier--Stokes equations. Hereby, the time discretization is chosen in such a way that the thermodynamical consistency of the system (cf. \eqref{p13energyIE}) is maintained.

We ensure the existence of feasible and globally optimal points for the respective optimal control problem and provide a first characterization of those points via a stationarity system of limiting $\mathcal{E}$-almost C-stationary type.
We proceed with a thorough analysis of the sensitivity and differentiability properties of the associated control-to-state operator which culminates in the presentation of a strong stationarity system.

Our analytical results are subsequently supplemented by the development and demonstration of two numerical solution algorithms, which compute discrete approximations of C-stationary or strong stationary points of the optimal control problem \eqref{p13optprob} below.
In order to handle the tremendous computational effort caused by repeatedly solving the large scale Navier-Stokes systems, we incorporate an adaptive mesh refinement strategy based on a goal-oriented error estimator.

\subsection{The semi-discrete CHNS system and the optimal control problem}
\label{p13sec:OC}

Let us start by presenting the underlying time discretization of the CHNS system and imposing some common assumptions on the related physical data.
For this purpose, we choose an arbitrary time step-size $\tau>0$ and denote the total number of time  {instants} by $K\in\mathbb{N}$ .
Moreover, we introduce a distributed force $u$ on the right-hand side of the Navier--Stokes equations.
\begin{defn}[Semi-discrete CHNS system]\label{p13defsemidis}
For a given initial state $(\varphi _{-1}, v _{0})=(\varphi _a , v _a) \in \left({H}^2_{\partial_n}(\Omega)\cap \mathbb{K}\right) \times  H^{2}_{0,\sigma}(\Omega;\mathbb{R}^n)$ we say that a triple
\begin{align*}
(\varphi,\mu,v)=((\varphi_i)_{i=0}^{K-1},(\mu_i)_{i=0}^{K-1},(v_i)_{i=1}^{K-1})
\end{align*}
in ${H}^2_{\partial_n}(\Omega)^{K} \times {H}^2_{\partial_n}(\Omega)^{K}\times H^{1}_{0,\sigma}(\Omega;\mathbb{R}^n)^{K-1}$
solves the semi-discrete CHNS system
with respect to a given control
$u=(u_i)_{i=1}^{K-1}\in L^2(\Omega;\mathbb{R}^n)^{K-1}$,
if it holds for all $\phi\in \overline{H}^1(\Omega)$ and $\psi\in H^{1}_{0,\sigma}(\Omega;\mathbb{R}^n)$ that
\begin{subequations}
\begin{align} %
&\left\langle\frac{\varphi_{i+1} -\varphi_{i} }{\tau},\phi\right\rangle
+\left\langle v_{i+1}\nabla\varphi_{i},\phi\right\rangle
 {+\left( m(\varphi_{i})\nabla\mu_{i+1},\nabla\phi\right)}=0,\label{p13firsttim1}\\ 
&\hspace*{0.6cm}  {\left( {\nabla}\varphi_{i+1},\nabla\phi\right)}
+\left\langle a_{i+1},\phi\right\rangle
-\left\langle \mu_{i+1},\phi\right\rangle
-\left\langle \kappa\varphi_{i},\phi\right\rangle= 0,\label{p13firsttim2}\\ 
&\left\langle\frac{\rho(\varphi_{i}) v_{i+1}-\rho(\varphi_{i-1}) v_i}{\tau},\psi\right\rangle_{H^{-1}_{0,\sigma},H^1_{0,\sigma}}
 {-\left( v_{i+1}\otimes \rho(\varphi_{i-1})v_i,\nabla\psi\right)}\nonumber\\
& {+\left( v_{i+1}\otimes \frac{\rho_2-\rho_1}{2}m(\varphi_{i-1})\nabla\mu_i,\nabla\psi\right)}
+(2\eta(\varphi_{i})\epsilon(v_{i+1}),\epsilon (\psi))\nonumber\\ %
&\hspace*{4cm}-\left\langle\mu_{i+1}\nabla\varphi_{i},\psi\right\rangle_{H^{-1}_{0,\sigma},H^1_{0,\sigma}}
=\left\langle u_{i+1},\psi\right\rangle_{H^{-1}_{0,\sigma},H^1_{0,\sigma}},\label{p13firsttim3} 
\end{align}\label{p13firsttim}
\end{subequations}
\hspace{-2ex} with $a_i\in \partial\Psi_0(\varphi_{i})$.
The first two equations are supposed to hold for every $0\leq i+1 \leq K-1$ and 
the last equation holds for every $1\leq i+1 \leq K-1$. 

The corresponding solution operator is denoted by $S_\Psi$, i.e. $(\varphi,\mu,v)\in S_\Psi(u)$. 
\end{defn}
In the above definition, the boundary conditions \eqref{p13eq:bdryCH} and the solenoidality of the velocity field \eqref{p13eq:NaSt2} are integrated in the chosen function spaces
\begin{align*}
H^{k}_{0,\sigma}(\Omega;\mathbb{R}^n)&:=\left\{f\in H^{k}(\Omega;\mathbb{R}^n)\cap H^{1}_0(\Omega;\mathbb{R}^n):\textnormal{div} f =0,\text{ a.e. on }\Omega\right\},\\ 
{H}^{k}_{\partial_n}(\Omega) &:= 
\left\{f\in {H}^{k}(\Omega):\partial_n f_{|\partial\Omega}=0\text{ on }{\partial\Omega}\right\},\ {k\geq 2}, 
\end{align*} 
for $\varphi,\mu$ and $v$.
Furthermore, the definition already includes the inherent regularity properties of $\varphi$ and $\mu$ which anticipates the results of Theorem \ref{p13exsols1} below.

Moreover, the semi-discrete CHNS system involves three time  {instants} $(i-1,i,i+1)$ and $(\varphi_0,\mu_0)$ is characterized in an initialization step by the (decoupled) Cahn--Hilliard system only. At the subsequent time  {instants}, the strong coupling of the Cahn--Hilliard and Navier--Stokes system is maintained.

In this work, we consider non-degenerate mobility and viscosity coefficients $m,\eta\in C^2(\mathbb{R})$, i.e. $0<c_1\leq \min_{x\in\mathbb{R}}\{m(x),\eta(x)\}$.
We further assume that $m$ and $\eta$, as well as their derivatives up to second order are bounded, which is typically satisfied if they originate from a practical application.

As noted above, the free energy density is related to the double-obstacle potential.
In other words,
the functional
$\Psi_0:{H}^1(\Omega)\rightarrow\mathbb{R}$ 
is given by $\Psi_0(\varphi):=\int_\Omega \iota_{[\psi_1;\psi_2]}(\varphi(x)) dx$,
where $\iota_{[\psi_1;\psi_2]}$ denotes the indicator function of ${[\psi_1;\psi_2]}$, i.e.
\begin{align}
\iota_{[\psi_1;\psi_2]}:=\left\{\begin{array}[c]{ll}
+\infty & \text{if } z< \psi_1,\\
0& \text{if } \psi_1\leq z\leq\psi_2, \\
+\infty & \text{if } z>\psi_2,\\
\end{array} \right. \quad \psi_1<0<\psi_2.\label{p13dob} 
\end{align}

As a consequence, the inclusion \eqref{p13firsttim2} ensures that the order parameter $\varphi_i$ is contained in ${[\psi_1;\psi_2]}$ almost everywhere (a.e.) on $\Omega$ for every time instance $-1\leq i \leq K-1$ assuming that
the initial data is well-posed in the sense that
\begin{align}
\varphi_a\in\mathbb{K}:=\left\{v\in {H}^1(\Omega):\psi_1\leq v\leq\psi_2\textnormal{ a.e. in }\Omega\right\}. 
\end{align}

In order to formulate the associated optimal control problem to \eqref{p13firsttim}, we introduce an objective functional $\mathcal{J}:\mathcal{X}\rightarrow\mathbb{R}$ defined on
$$\mathcal{X}:= {H}^1(\Omega)^{K} \times {H}^1(\Omega)^{K}\times H^{1}_{0,\sigma}(\Omega;\mathbb{R}^n)^{K-1}\times L^2(\Omega;\mathbb{R}^n)^{K-1},$$
and assume that $\mathcal{J}$ is convex, weakly lower-semi-continuous, Fr{\'e}chet differentiable, and partially coercive.
\begin{defn}\label{p13optprob}
We study the optimal control problem
\begin{align}
\begin{aligned}
&\min \mathcal{J}(\varphi,\mu,v,u)\textnormal{ over } (\varphi,\mu,v,u)\in\mathcal{X}\\
&\textnormal{s.t. }(\varphi,\mu,v)\in S_\Psi(u). 
\end{aligned} \label{p13optprob.Ppsi} 
\end{align}
\end{defn}
For our numerical computations below, we consider the specific functional
 \begin{align}
  \mathcal{J}(\varphi,\mu,v,u):=\frac{1}{2}\left\|\varphi_{K-1}-\varphi_{d} \right\|^2+\frac{\xi}{2}\left\| u \right\|^2,\ \xi>0,\label{p13trackobj}
 \end{align}
where $\varphi_{d}\in L^2(\Omega)$ represents a desired state.
The so-called tracking type functional, which is used in various applications, clearly satisfies the above assumptions.

\subsection{Existence of feasible and globally optimal points}
\label{p13sec:exist}

One of the main requirements for the existence of solutions to \eqref{p13optprob.Ppsi} is the boundedness of the state.
In our setting, this property follows from the energetic stability of the chosen discretization in time.
More precisely, we have the following (dissipative) energy law for the total energy
\begin{align}
 E(v,\varphi,\varphi_{-1})=\int_\Omega \rho(\varphi_{-1})\frac{|v|^2}{2} dx + \int_\Omega \frac{|\nabla\varphi|^2}{2} dx + \Psi(\varphi).
\end{align}
associated with the semi-discrete CHNS system \eqref{p13firsttim}.
\begin{lem}[Energy estimate for a single time step]\label{p13energyest}
Let ${\varphi_i},\varphi_{i-1}\in {H}^2_{\partial_n}(\Omega) \cap \mathbb{K}$, $\mu_i\in {H}^2_{\partial_n}(\Omega)$,
${v_i}\in H^{1}_{0,\sigma}(\Omega;\mathbb{R}^n)$
and $u_{i+1}\in (H^{1}_{0,\sigma}(\Omega;\mathbb{R}^n))^*$ be given. 

If $(\varphi_{i+1},\mu_{i+1},v_{i+1})\in {H}^1(\Omega) \times {H}^1(\Omega)\times H^{1}_{0,\sigma}(\Omega;\mathbb{R}^n)$
satisfies the system \eqref{p13firsttim},
then the corresponding total energy is bounded by
\begin{align}
&E(v_{i+1},\varphi_{i+1},\varphi_{i})
+\int_\Omega \rho(\varphi_{i-1}) \frac{\left| v_{i+1} -v_i \right|^2}{2}dx
+\int_\Omega \frac{\left| \nabla\varphi_{i+1} -\nabla\varphi_i \right|^2}{2}dx\nonumber\\
&\hspace{1cm}+\tau\int_\Omega2\eta(\varphi_i)\left|\epsilon(v_{i+1})\right|^2dx
+\tau\int_\Omega m(\varphi_i)\left|\nabla\mu_{i+1}\right|^2dx
+\int_\Omega \kappa\frac{(\varphi_{i+1}-\varphi_i)^2}{2}\nonumber\\
&\hspace{2cm}\leq E(v_{i},\varphi_{i},\varphi_{i-1})
+\left\langle u_{i+1},v_{i+1}\right\rangle_{H^{-1}_{0,\sigma},H^1_{0,\sigma}}.\label{p13EE}
\end{align}
\end{lem}
It should be noted that the density is always positive, since $\varphi_i$ is contained in $\mathbb{K}$ for every $i$.
Consequently, all the terms of the left-hand side of the inequality are always non-negative 
such that Lemma \ref{p13energyest} indeed ensures that the energy of the next time step is non-increasing if the external force $u_{i+1}$ is absent.

Lemma \ref{p13energyest} allows us to verify the existence of solutions to the CHNS system \eqref{p13firsttim} via the repeated application of Schaefer's fixed point theorem.
The proof further involves arguments from PDE theory and monotone operator theory.

\begin{thm}[Existence of feasible points]\label{p13exsols1}
Let $u\in L^2(\Omega;\mathbb{R}^n)^{K-1}$ be given.

Then the semi-discrete CHNS system admits a solution
$(\varphi,\mu,v)\in {H}^2_{\partial_n}(\Omega)^{K} \times {H}^2_{\partial_n}(\Omega)^{K}\times H^{2}_{0,\sigma}(\Omega;\mathbb{R}^n)^{K-1}$.
\end{thm}

The last theorem also ensures an additional regularity of the state, which is necessary to guarantee that the system \eqref{p13firsttim} is well-posed for each time step.
The proof relies on the regularity theory for Navier-Stokes equations and variational inequalities.

By Theorem \ref{p13exsols1} the feasible set of problem \eqref{p13optprob.Ppsi} is non-empty.
Then the existence of globally optimal points can be verified via standard arguments from optimization theory.
\begin{thm}[Existence of global solutions]\label{p13exsol}
The optimization problem \eqref{p13optprob.Ppsi} possesses a global solution. 
\end{thm}
For more details on the results presented in this subsection, we refer the reader to \cite{Hintermueller2015}.

\subsection{$\mathcal{E}$-almost C-stationary points}
\label{p13sec:Cstat}

After securing the existence of solutions to the optimal control problem \eqref{p13optprob.Ppsi} we target a more precise characterization of globally and/or locally optimal points via necessary optimality conditions.
This lays the foundation to the development of efficient numerical solution methods in the subsequent subsections.

As a first step, we establish a limiting $\mathcal{E}$-almost C-stationarity system.
For this purpose, we additionally assume that $\mathcal{J}'$ is a bounded mapping
		and $\frac{\partial \mathcal{J}}{\partial u}$ satisfies the following weak lower-semicontinuity property
		\vspace{-0.1cm}
	\[
		\Big\langle \frac{\partial \mathcal{J}}{\partial u}(\hat z), \hat u \Big\rangle
		\>\le\>	\liminf_{k\rightarrow\infty}
				\Big\langle \frac{\partial \mathcal{J}}{\partial u}(\hat z^{(k)}), \hat u^{(k)} \Big\rangle {,}
	\]
	where $\hat z^{(k)}$ converges weakly in 
		${H}^2_{\partial_n }(\Omega )^ K \times {H}^2_{\partial_n }(\Omega )^ K \times H^1_{0,\sigma}(\Omega ;\mathbb{R}^n )^{K-1} \times L^2(\Omega ;\mathbb{R}^n )^{K-1}$
		towards a limit point $\hat z$. 
	Here and in the following $z$ represents the primal variables, i.e. $\hat z:=(\hat\varphi , \hat\mu , \hat v , \hat u )$.

The derivation is based on a penalization of the lower-level problem, where the double-obstacle potential is approximated by certain smooth double-well type potentials $\Psi_k,\ k\in\mathbb{N}$.
This gives rise to a family of smooth auxiliary nonlinear programs ($P_{\Psi_k}$) for which the following necessary optimality system can be derived via a well-known result from Zowe and Kurcyusz \cite[Theorem 4.1]{Zowe1979}.

\begin{thm}[First-order optimality conditions for smooth potentials]\label{p13T:Mult}
	Let
		$\overline {z}$ be a minimizer of the auxiliary problem ($P_{\Psi_k}$).
		
	Then there exist $( p,r,q,\lambda )\in {H}^1(\Omega )^ K \times {H}^1(\Omega )^{K} \times { H^1_{0,\sigma}(\Omega ;\mathbb{R}^n ) }^{K-1}\times {\overline{H}^1(\Omega )^*}^{K} $, with
		$\lambda _{i}:=\Psi_k'' (\varphi _{i+1} )^* r _{i}$,
		such that
	\begin{align}	
							- \frac1 \tau ( p _{i} - p _{i-1} ) + m'(\varphi_{i})\nabla \mu _{i+1} \cdot\nabla p _{i} - \mathop{\rm div}( p _{i} v _{i+1} ) - \Delta r _{i-1} 				\nonumber \\ \hskip1cm
								+ \lambda _{i-1}
								- \kappa r _{i+1} - \frac 1 \tau \rho' (\varphi_{i}) v _{i+1} \cdot( q _{i+1} - q _{i} )	 							\nonumber \\ \hskip1cm
								- (\rho' (\varphi_{i}) v _{i+1}  -\frac{\rho_2-\rho_1}{2}m'(\varphi_{i})\nabla \mu _{i+1} ) (D q _{i+1} )^\top v _{i+2} 									\nonumber \\ \hskip1cm
								+ 2\eta' (\varphi_{i}) \epsilon( v _{i+1} ) : D q _{i} + \mathop{\rm div}( \mu _{i+1} q _{i} )									 
					& \>=\>		\frac {\partial \mathcal{J}}{\partial \varphi _{i} }(\overline {z}),														\label{p13T:Mult.1}\\ 
							- r _{i-1} -\mathop{\rm div}( m(\varphi_{i-1}) \nabla p _{i-1} )													
								- \mathop{\rm div}( \frac{\rho_2-\rho_1}{2}m(\varphi_{i-1}) (D q _{i} )^\top v _{i+1} )												\nonumber \\ \hskip1cm
								- q _{i-1} \cdot\nabla \varphi _{i-1} 											 
					& \>=\>	\frac {\partial \mathcal{J}}{\partial \mu _{i} }(\overline {z}),	\label{p13T:Mult.2} 
  \end{align}

 \begin{align}
 							- \frac 1 \tau \rho (\varphi_{j-1}) ( q _{j} - q _{j-1} ) - \rho (\varphi_{j-1}) (D q _{j} )^\top v _{j+1} 									\nonumber \\ \hskip1cm
								- (D q _{j-1} )( \rho (\varphi_{j-2}) v _{j-1}  -\frac{\rho_2-\rho_1}{2}m(\varphi_{j-2}) \nabla \mu _{j-1} )									\nonumber \\ \hskip1cm
								- \mathop{\rm div}( 2 \eta (\varphi_{j-1}) \epsilon( q _{j-1} ) ) + p _{j-1} \nabla \varphi _{j-1} 								 
					& \>=\>	\frac {\partial \mathcal{J}}{\partial v _{j} }(\overline {z}),															\label{p13T:Mult.3}\\ 
 %
							\frac {\partial \mathcal{J}}{\partial u _{j} }(\overline {z}) -  q _{j-1} 
					& \>=\>	0 \label{p13T:Mult.4}
	\end{align}
	for all $i=0,...,{K-1}$ and $j=1,...,{K-1}$.
 %
	Here, we use the convention that $ p _{i} , r _{i} , q _{i} $ are equal to $0$ for $i\ge{K-1}$ along with $ q _{-1}$ and
		$\varphi _{i} ,\mu _{i} , v _{i} $ for $i\ge K $.
\end{thm}
A careful limit analysis with respect to a vanishing penalization parameter yields the following stationarity system for the optimal control problem \eqref{p13optprob.Ppsi}, cf. \cite{Hintermueller2015}.

\begin{thm}[Limiting $\mathcal{E}$-almost C-stationarity]\label{p13T:Doub}
 Let $( \varphi ^{(k)} , \mu ^{(k)} , v ^{(k)} , u ^{(k)} )$
 be a minimizer for (P$_{ \Psi_{k} }$)
		and let further $( p ^{(k)} , r ^{(k)} , q ^{(k)}, \lambda^{(k)} )$ 
		be given as in Theorem~\ref{p13T:Mult}.

	Then there exists a weakly convergent subsequence 
	\begin{align}
	&\left\{( \varphi ^{(m)} , \mu ^{(m)} , v ^{(m)} , u ^{(m)} , p ^{(m)} , r ^{(m)} , q ^{(m)},\lambda^{(m)} )\right\}_{m\in\mathbb{N}}\nonumber\\
	&\subset {H}^2_{\partial_n }(\Omega )^ K \times {H}^2_{\partial_n }(\Omega )^ K \times H^1_{0,\sigma}(\Omega ;\mathbb{R}^n )^{K-1} \times L^2(\Omega ;\mathbb{R}^n )^{K-1}\nonumber\\
	&\quad \times {{H}^1(\Omega )}^ K \times { {H}^1(\Omega )}^ K \times { H^1_{0,\sigma}(\Omega ;\mathbb{R}^n )}^{K-1}\times {{H}^1(\Omega )^*}^K
	\end{align}
	
	$ $
	and the limit point $( \varphi , \mu , v , u , p , r , q, \lambda )$ satisfies the adjoint system \eqref{p13T:Mult.1}-\eqref{p13T:Mult.4}, as well as
	\begin{align}
			(\, a _{i} , r _{i-1} \,)_{L^2}		&	=	0	,&
		\liminf (\, \lambda ^{(m)}_{i} , r ^{(m)}_{i-1} \,)_{L^2}	&	\ge	0	\label{p13T:Doub.2}.
	\end{align}
	Moreover,
		for every $\varepsilon >0$
		there exist a measurable subset $ M ^\varepsilon _{i} $ of $ M _{i} :=\{x\in\Omega \>:\> \psi_1 < \varphi _{i} (x)< \psi_2 \}$ with
		$| M _{i} \setminus M ^\varepsilon _{i} |<\varepsilon $ and
	\begin{align}
			\langle \lambda _{i} ,v\rangle=0	\enspace \enspace \enspace \enspace \forall v\in \overline{H}^1(\Omega ),\enspace v|_{\Omega \setminus M ^\varepsilon _{i} }=0.\label{p13epscom}
	\end{align}
\end{thm}

The above stationarity conditions correspond to a function space version of C-stationarity, see, e.g., \cite{Hintermuller2009,Scheel2000}.
The proof of the last condition \eqref{p13epscom} is based on the application of Egorov's theorem, cf. \cite{Barbu1984}, which motivated the notion of $\mathcal{E}$-almost C-stationarity.

\subsection{Strong stationarity}
\label{p13sec:Sstat}

Starting from the C-stationarity system of the previous section, it is possible to derive a
more restrictive stationarity system for the problem \eqref{p13optprob.Ppsi} employing the directional differentiability of the control-to-state operator $S_\Psi$.
In this subsection, we consider the control of the semi-discrete CHNS system for a single time step, i.e. $K=2$ and $\varphi_{-1}, \varphi_0, \mu_0, v_0$ are given. This corresponds to an instantaneous control problem.

First, we verify that the solution operator $S_\Psi$ of the semi-discrete CHNS system is Lipschitz continuous. 

\begin{thm}[Lipschitz continuity of $S_\Psi$]\label{p13LIPSS}
The mapping $S_\Psi:{H}^{-1}_{0,\sigma}(\Omega)\rightarrow {H}^1(\Omega) \times {H}^1(\Omega)\times H^1_{0,\sigma}(\Omega;\mathbb{R}^N)$  is Lipschitz continuous.
\end{thm}
The proof follows a similar line of argumentation as Lemma \ref{p13energyest}. 
An immediate consequence of the above theorem is that the solutions to the constraint system are uniquely determined by the control $u$.

Although solution operators of variational inequalities are generally not Fr{\'e}chet differentiable,
we can now compute the directional derivative of $S_\Psi$ via the following theorem.
\begin{thm}\label{p13directder}
The directional derivative of $S_\Psi$ at $\hat u\in {H}^{-1}_{0,\sigma}(\Omega)$ with $S_\Psi(\hat u)=(\hat \varphi,\hat \mu,\hat v)$ in direction $h\in {H}^{-1}_{0,\sigma}(\Omega)$
is 
the unique solution $(\chi,w,\zeta)\in H^1(\Omega)\times H^1(\Omega)\times H^1_{0,\sigma}(\Omega;\mathbb{R}^N)$
of the system
\begin{subequations}
\begin{align} %
\chi\in T_\mathbb{K}(\hat \varphi)\cap {{a} ^+}^\bot\cap { a ^-}^\bot,\\
\left\langle-\Delta \chi - w,v-\chi\right\rangle\geq 0,\ \forall v\in T_\mathbb{K}(\hat \varphi)\cap { a ^+}^\bot\cap { a ^-}^\bot,\label{p13dirderVI1}\\
\left\langle\frac{\chi  }{\tau},\phi\right\rangle
+\left\langle \zeta\nabla\varphi_{0},\phi\right\rangle
+\left( m(\varphi_{0})\nabla w,\nabla\phi\right)
=0,\label{p13dirderVI2}\\
\left\langle\frac{\rho(\varphi_{0}) \zeta}{\tau},\psi\right\rangle_{H^{-1}_{0,\sigma},H^1_{0,\sigma}}
-\left( \zeta\otimes \rho(\varphi_{-1})v_0,\nabla\psi\right)\nonumber\\
+\left( \zeta\otimes \frac{\rho_2-\rho_1}{2}m(\varphi_{-1})\nabla\mu_0,\nabla\psi\right)
+(2\eta(\varphi_{0})\epsilon(\zeta),\epsilon (\psi))\nonumber\\ %
\hspace*{4cm}-\left\langle w\nabla\varphi_{0},\psi\right\rangle_{H^{-1}_{0,\sigma},H^1_{0,\sigma}}
-\left\langle h,\psi\right\rangle_{H^{-1}_{0,\sigma},H^1_{0,\sigma}}=0.
\end{align}
\end{subequations}
Here, $T_\mathbb{K}(\hat \varphi)$ represents the tangent cone of $\mathbb{K}$ at $\hat \varphi$ and $ {a^{+/-}} ^\bot:=\{\phi\in H^1(\Omega):\left\langle \phi, a^{+/-} \right\rangle=0\}$ is the orthogonal space associated with
$a^+(x):=\max\{a(x),0\}$ and $a^-(x):=\min\{a(x),0\}$.  
\end{thm}

Note that $a^+$ and $a^-$ can be interpreted as the multipliers to the constraints $\varphi\leq 1$ and $\varphi\geq -1$ and
the convex constraint set $T_\mathbb{K}(\hat \varphi)\cap { a ^+}^\bot\cap { a ^-}^\bot$ associated to the variational inequality \eqref{p13dirderVI1}-\eqref{p13dirderVI2} is also called the critical cone, cf. \cite{Mordukhovich2006a}. 
The proof of Theorem \ref{p13directder} combines arguments from Jarusek et al. in \cite{Jarusek2003} and PDE theory. 

With the help of the directional derivative of $S_\Psi$, we derive strong stationarity conditions for \eqref{p13optprob.Ppsi} by evaluating the B-stationarity condition of the reduced optimization problem 
\begin{align}
\min_{u\in L^2(\Omega;\mathbb{R}^N)}\overline{\mathcal{J}}(u):=\mathcal{J}(S_\Psi(u),u) \label{p13reducopt} 
\end{align}
for suitable test directions.
\begin{thm}
If $\hat u$ is an optimal control of \eqref{p13optprob.Ppsi}, then there exists an adjoint state
$(p,r, q)\in {H}^1(\Omega)\times {H}^1(\Omega)\times H^1_{0,\sigma}(\Omega;\mathbb{R}^N)$ and $\lambda\in {H}^1(\Omega)^*$ such that for all $\phi\in H^1(\Omega)$ and $\psi\in H^1_{0,\sigma}(\Omega;\mathbb{R}^N)$
it holds that
  \begin{align}
 \left\langle D_\varphi \mathcal{J}[z_0]+\frac{r  }{\tau},\phi\right\rangle+\left( \nabla p,\nabla\phi\right)+\left\langle\lambda,\phi\right\rangle&=0,\label{p13statsys1}\\
\left( m(\varphi_{0})\nabla r,\nabla\phi\right)-\left\langle p,\phi\right\rangle-\left\langle  q\nabla\varphi_{0},\phi\right\rangle&=0,\\
 \left\langle\frac{\rho(\varphi_{0})}{\tau} q,\psi\right\rangle_{H^{-1}_{0,\sigma},H^1_{0,\sigma}} 
 - \left\langle \nabla q\nu,\psi\right\rangle_{H^{-1}_{0,\sigma},H^1_{0,\sigma}}&\nonumber\\
 +\left\langle 2\eta(\varphi_{0})\epsilon( q),\epsilon(\psi)\right\rangle_{H^{-1}_{0,\sigma},H^1_{0,\sigma}}
 - \left\langle r\nabla\varphi_{0},\psi\right\rangle_{H^{-1}_{0,\sigma},H^1_{0,\sigma}}
 &=0,\\
 \left\langle- q,\psi\right\rangle_{H^{-1}_{0,\sigma},H^1_{0,\sigma}}
 + \left\langle D_u \mathcal{J}[\hat{z}],\psi\right\rangle_{H^{-1}_{0,\sigma},H^1_{0,\sigma}}
 &=0,\label{p13statsys4}\\
 \lambda\in \left(T_\mathbb{K}(\hat \varphi)\cap { a ^+}^\bot\cap { a ^-}^\bot\right)^0,\label{p13statsys5}\\
  q
\in
\left(\left[D\left(\left(T_\mathbb{K}(\hat \varphi)\cap { a ^+}^\bot\cap { a ^-}^\bot\right)^0\times H^1_{0,\sigma}(\Omega;\mathbb{R}^N)\right)\right]_2\right)^0,\label{p13statsys6}
\end{align}
where $D$ is a specific linear operator and the subscript $K^0$ represents the polar cone of the cone $K$.
\end{thm}
This concludes our analytical investigations.
We point out that the strong stationarity conditions represent the most selective stationarity system available for the problem under consideration up to this point in time.

\subsection{Adaptive mesh refinement}
\label{p13sec:adapt}

In the following subsections, we discuss efficient numerical solution methods for the problem \eqref{p13optprob.Ppsi},
where the objective functional is given by \eqref{p13trackobj}, based on our analytical results.
The main challenges hereby are imposed by the non-differentiability of the solution operator due to the Cahn-Hilliard system and the immense numerical expense caused by repeatedly solving the large scale Navier-Stokes type primal and dual systems.

We deal with the second challenge by developing a goal-oriented error estimator based on the dual-weighted residual approach, cf., e.g., \cite{Brett2015}.
This allows us to implement an adaptive mesh refinement strategy, which acknowledges the error contributions of the primal residuals, the dual residuals and the mismatch in the complementarity terms, to reduce the computational effort. 

The central idea of this approach is depicted by the subsequent theorem, which
estimates the difference of the objective values at stationary points of the semi-discrete and the fully discretized problem with the help of the associated MPCC-Lagrangian $\mathcal{L}$, cf. \cite{HHKK}.

\begin{thm}\label{p13errtheo}
Let 
$(y,u,\Phi,\pi,\lambda^{+},\lambda^{-})$ be a stationary point of the optimal control problem \eqref{p13optprob.Ppsi} 
and assume that $(y_h,u_h,\Phi_h,\pi_h,\lambda_{h}^+,\lambda_{h}^-)\in\mathcal{Y}_h$ 
satisfies the discretized stationarity system.
Then it holds that
\begin{align}
  \mathcal{J}(\varphi_h,\mu_h,v_h,u_h)
  - \mathcal{J}(\varphi,\mu,v,u)
  = \frac{1}{2}\left(\sum_{i=0}^{K-1}\left\langle a_h^{i},\pi^i\right\rangle
  -\sum_{i=0}^{K-1}\left\langle a^{i},\pi_h^{i}\right\rangle\right)\nonumber\\
  -\frac{1}{2}\left(\sum_{i=0}^{K-1}\left\langle(\lambda^i)^+,\varphi_h^{i}-\psi_2\right\rangle
  -\sum_{i=0}^{K-1}\left\langle (\lambda_{h}^i)^+,\varphi^i-\psi_2\right\rangle\right)\nonumber\\
  +\frac{1}{2}\left(\sum_{i=0}^{K-1}\left\langle (\lambda^i)^-,\varphi^i_h-\psi_1\right\rangle
  -\sum_{i=0}^{K-1}\left\langle (\lambda_{h}^i)^-,\varphi^i-\psi_1\right\rangle\right)\nonumber\\
  +\frac{1}{2}\nabla_x \mathcal{L}(y_h,u_h,\Phi_h,\pi_h,\lambda_{h}^+,\lambda_{h}^-) ((y_h,u_h,\Phi_h)-(y,u,\Phi))\nonumber\\
   {+O\left(\|(y_h,u_h,\Phi_h)-(y,u,\Phi)\|^3\right)},\label{p13firstest}
\end{align}
where $O$ denotes the Landau symbol Big-$O$. 
\end{thm}
This allows us to approximate the discretization error with respect to the objective function as follows
\begin{equation}
  \label{p13eq:DWR:finalErrRep}
  \begin{aligned}
    \mathcal{J}(\varphi_h,\mu_h,v_h,u_h)
    -& \mathcal{J}(\varphi,\mu,v,u)\\
    \approx\sum_{i=0}^{K-1}&(\eta_{CM1,i}+\eta_{CM2,i}+\eta_{CM3,i}+\eta_{CM4,i}
    +\eta_{CH1,i}\\
    &+\eta_{CH2,i}+\eta_{NS,i}+\eta_{AD\varphi,i}+\eta_{AD\mu,i}+\eta_{ADv,i}),
  \end{aligned}
\end{equation}
where the complementarity error terms $\eta_{CM1,i},..,\eta_{CM4,i}$, the weighted primal residuals $\eta_{CH1,i},\eta_{CH2,i},\eta_{NS,i}$ and the weighted dual residuals $\eta_{AD\varphi,i},\eta_{AD\mu,i},\eta_{ADv,i}$
are defined as in \cite[Section 4]{HHKK}.
These individual error terms can be evaluated separately on each patch of the current mesh due to their integral structure.
In order to obtain a fully a-posteriori error estimator the continuous quantities are approximated with the help of a
local higher-order approximation based on the respective discrete variables. 

\subsection{Penalization algorithm}
\label{p13sec:penalgo}

A first approach to handle the non-differentiability of $S_\Psi$ numerically is motivated by the penalization method of Subsection \ref{p13sec:Cstat}.
Namely, we solve a sequence of auxiliary optimization problems, 
where we approximate $\Psi_0$ by 
\begin{align*}
  \Psi_{0,\alpha}(\varphi) := \frac{{1}}{2\alpha}\left(
  \max(0,\varphi-1)^2 + \min(\varphi+1)^2
  \right),\ \alpha>0,\ \alpha \rightarrow 0.
\end{align*}
The resulting nonlinear programs can be solved by a standard
steepest descent method 
and the calculated solution approximates a C-stationary point of \eqref{p13optprob.Ppsi}, if the complementarity conditions of Theorem \ref{p13T:Doub} are satisfied sufficiently well, i.e. up to a given tolerance $tol_c$.
In combination with an outer adaptation loop based on the error estimator \eqref{p13eq:DWR:finalErrRep}, this yields the Algorithm~\ref{p13n15:alg:alg:overallAlgorithm}.

\newpage

\begin{algorithm}[H]
\KwData{Initial data: $\varphi_{a},v_a$;} 
\Repeat{$\sum_{i=1}^K |\mathcal T^i| > \mathcal{A}_{\max}$} 
{
\Repeat{\textnormal{complementarity conditions are satisfied up to a tolerance $\epsilon_{tol}$}}
{
solve the regularized problem $(P_{\Psi_\alpha})$ using a steepest descent method; \label{p13alg:alg:solveP} 

decrease $\alpha$;
}
\label{p13alg:alg:optLoop_end} 
calculate the error indicators and identify the sets $\mathcal M_r${,}$\mathcal M_c$ of cells to
refine/coarsen\; \label{p13alg:alg:EstAdap}
adapt $(\mathcal T^i)_{i=1}^K$ based on $\mathcal M_r$ and $\mathcal M_c$\;
}\label{p13alg:alg:finishType}
 \caption{The overall solution procedure}
\label{p13n15:alg:alg:overallAlgorithm}
\end{algorithm}

Hereby, the outer adaptation loop
relies on the D\"orfler marking procedure.
Hence, the error indicators from \eqref{p13eq:DWR:finalErrRep}
are evaluated 
for all time steps $i$ and for all cells $T\in \mathcal{T}^i$ of the current triangulation $(\mathcal T^i)_{i=1}^K$.
Then we choose a set $\mathcal M_r$ of cells to be refined as the set with the smallest cardinality which satisfies
\begin{align*}
  \sum_{T\in \mathcal M_r}  \eta_T \geq \theta^r \sum_{i=1}^K\sum_{T\in\mathcal
  T^i}
  \eta_T,
\end{align*}
for a given parameter $0<\theta^r < 1$. 
Due to the movement of the interface, we also select cells for coarsening
if the
calculated error indicator is smaller than a certain fraction of the mean error, i.e.
\begin{align*}
  M_c := \left\{ T \in (\mathcal T^i)_{i=1}^K\,|\, 
  \eta_T \leq \frac{\theta^c}{\mathcal{A}} \sum_{i=1}^K\sum_{T\in\mathcal T^i} \eta_T
  \right\},
\end{align*}
where $0<\theta^c<1$ is fixed and $\mathcal{A} := \sum_{i=1}^K |\mathcal T^i|$.
The mesh refinement process is terminated if a desired total number of cells $\mathcal{A}_{\max}$ is exceeded. 

Moreover, the problem is discretized in space using Taylor-Hood finite elements,
i.e. we utilize linear finite elements for
$\varphi$, $\mu$, and $p$ and quadratic finite elements for $v$.
For more details on the implementation of the algorithm and the numerical results we refer to \cite{HHKK}.

\begin{figure}
  
  \centering
  
  \fbox{
  \includegraphics[width=0.2\textwidth]{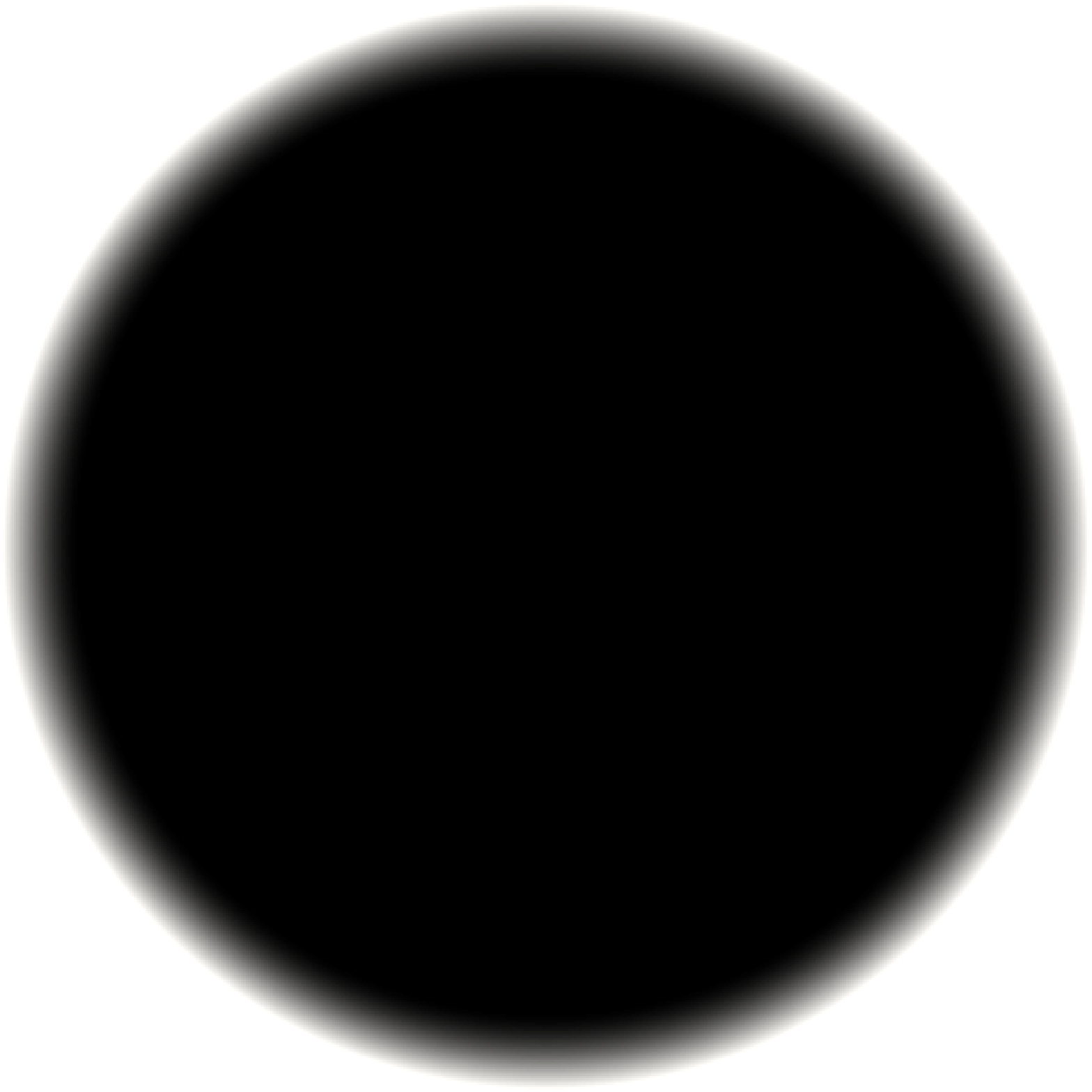}
  }
  \hspace{0.5cm}
  \fbox{
  \includegraphics[width=0.2\textwidth]{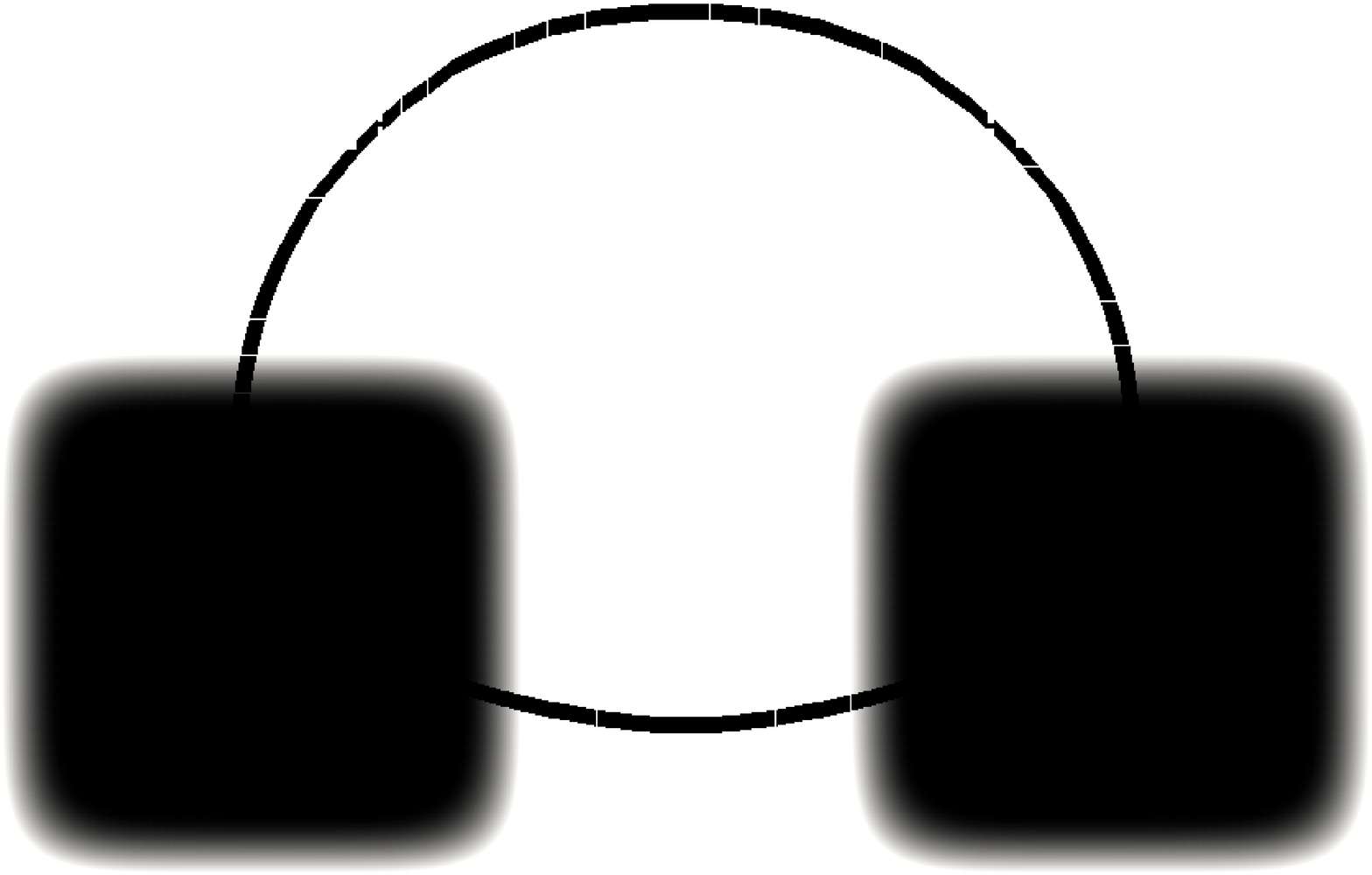}
  }
  \hspace{0.5cm}
  \fbox{
  \includegraphics[width=0.2\textwidth]{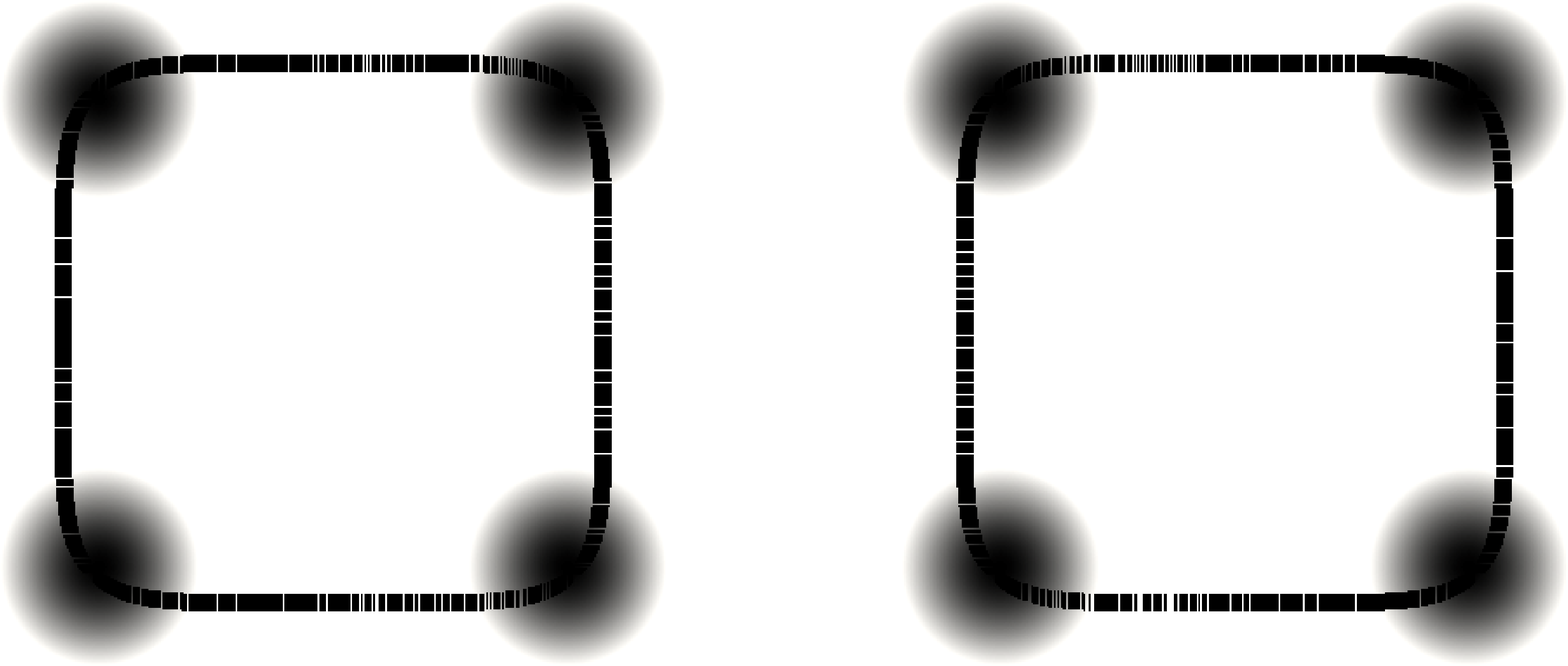}
  }

  \caption{The initial shape $\varphi_0$, the desired shape $\varphi_d$, the ansatz for the control $u$.}
  \label{p13fig:num:c2sX:c0cdBu}
\end{figure}

Let us briefly illustrate the performance of the proposed Algorithm \ref{p13n15:alg:alg:overallAlgorithm} at the hands of a specific example.
Our goal is to control the motion of a circular bubble to prevent it from rising and split it into two
square shaped bubbles. 
For this purpose, $2\times 4$ locally supported Ansatz functions of the control are distributed over the two-dimensional domain as depicted in Figure \ref{p13fig:num:c2sX:c0cdBu}.
The figure further shows the initial state $\varphi_a$, the desired shape $\varphi_d$
together with the zero level line of the phase field at final time if no
control is applied.
The corresponding objective functional is defined as in \eqref{p13trackobj} with $\xi = 1e-11$.

The associated fluid parameters are given by $\rho_1 = 1000$, $\rho_2 = 100$, $\eta_1 = 10$, $\eta_2 = 1$, and $\sigma = 24.5 \cdot\frac{2}{\pi}$ and
are taken from a benchmark problem for rising bubble dynamics in
\cite{n15_2009_HysingEtAl}. Furthermore, we incorporate a gravitational acceleration $g=0.981$ in the vertical direction and set $\epsilon = 0.02$,
$m(\varphi) \equiv \frac{1}{25000}$.
The time horizon is set to $T=1.0$ and the time step size is $\tau = 0.00125$.

For the marking procedure we use the parameters 
$\theta^r = 0.7$ and $\theta^c = 0.01$.
Furthermore, the stopping criteria use the tolerance
$tol_c = 1e-3$ for the complementarity conditions and the maximum amount of cells $\mathcal{A}_{max} = 8e6$ for the adaptation process, which relates to $1e4$ cells in average per time instance.

The optimal solution on the first level and for the initial value for $\alpha$ is found after 26
steepest descent iterations, while the complete algorithm terminates after 419 steepest descent
steps.
Hereby, the algorithm solves the auxiliary optimization problems
10 times, i.e. line \ref{p13alg:alg:solveP} of
Algorithm~\ref{p13n15:alg:alg:overallAlgorithm}
is executed 10 times.
After the first two solves the Moreau--Yosida parameter was decreased,
and after 
the next 8 solves the algorithm directly proceeded with outer adaptation loop.

\begin{figure}
  \centering
  \hfill
  \fbox
  {
  \includegraphics[width=0.15\textwidth]{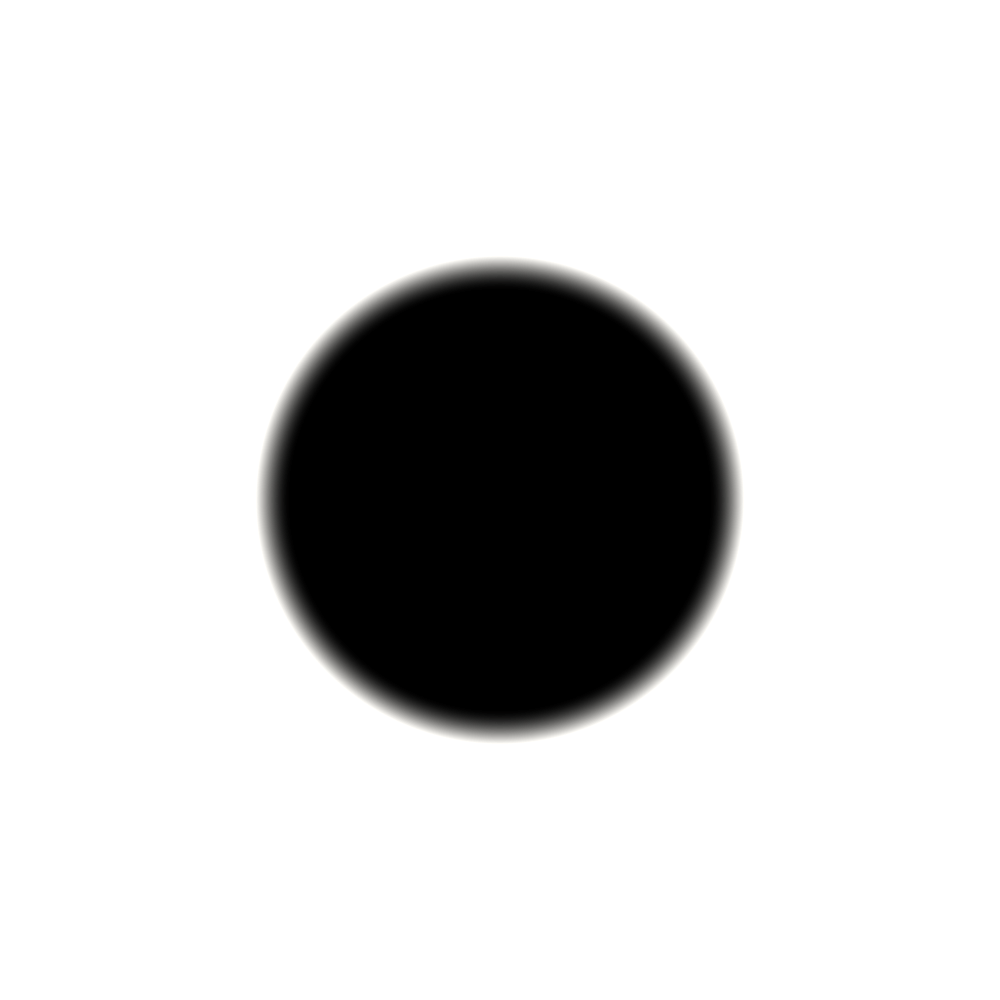}
  }
  \hfill
  \fbox
  {
  \includegraphics[width=0.15\textwidth]{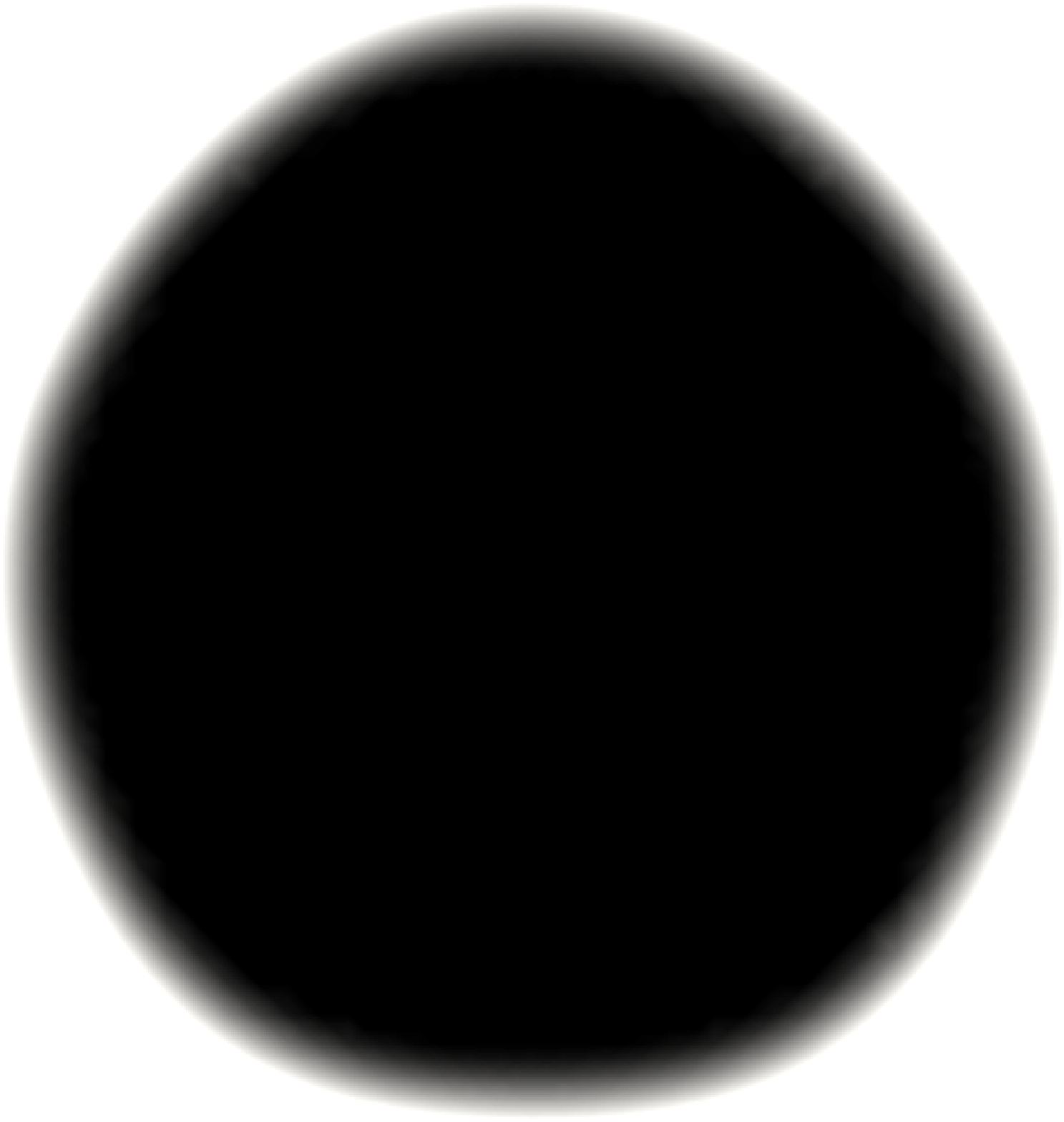}
  }
  \hfill
  \fbox
  {
  \includegraphics[width=0.15\textwidth]{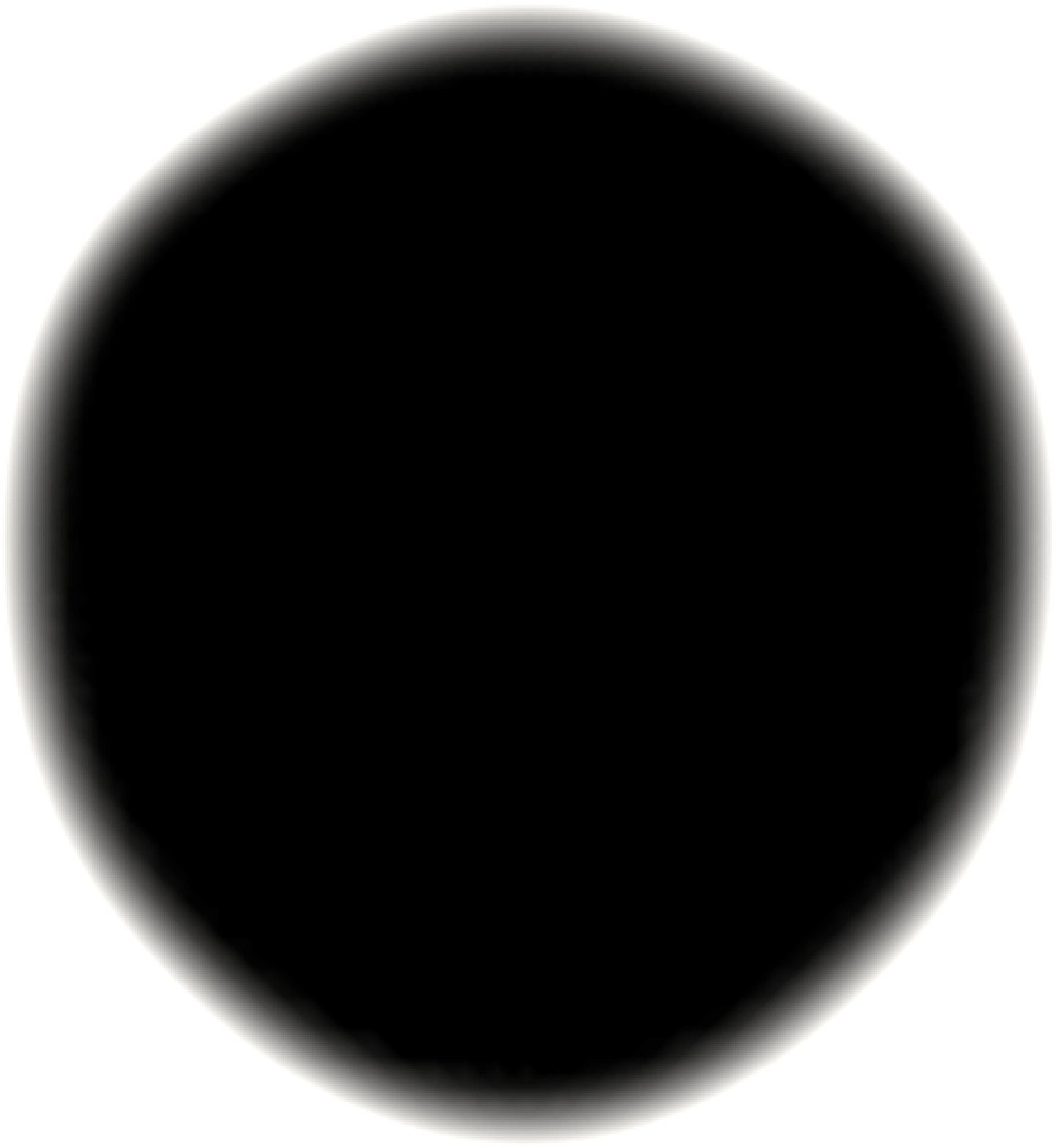}
  }
  \hfill
  \fbox
  {
  \includegraphics[width=0.15\textwidth]{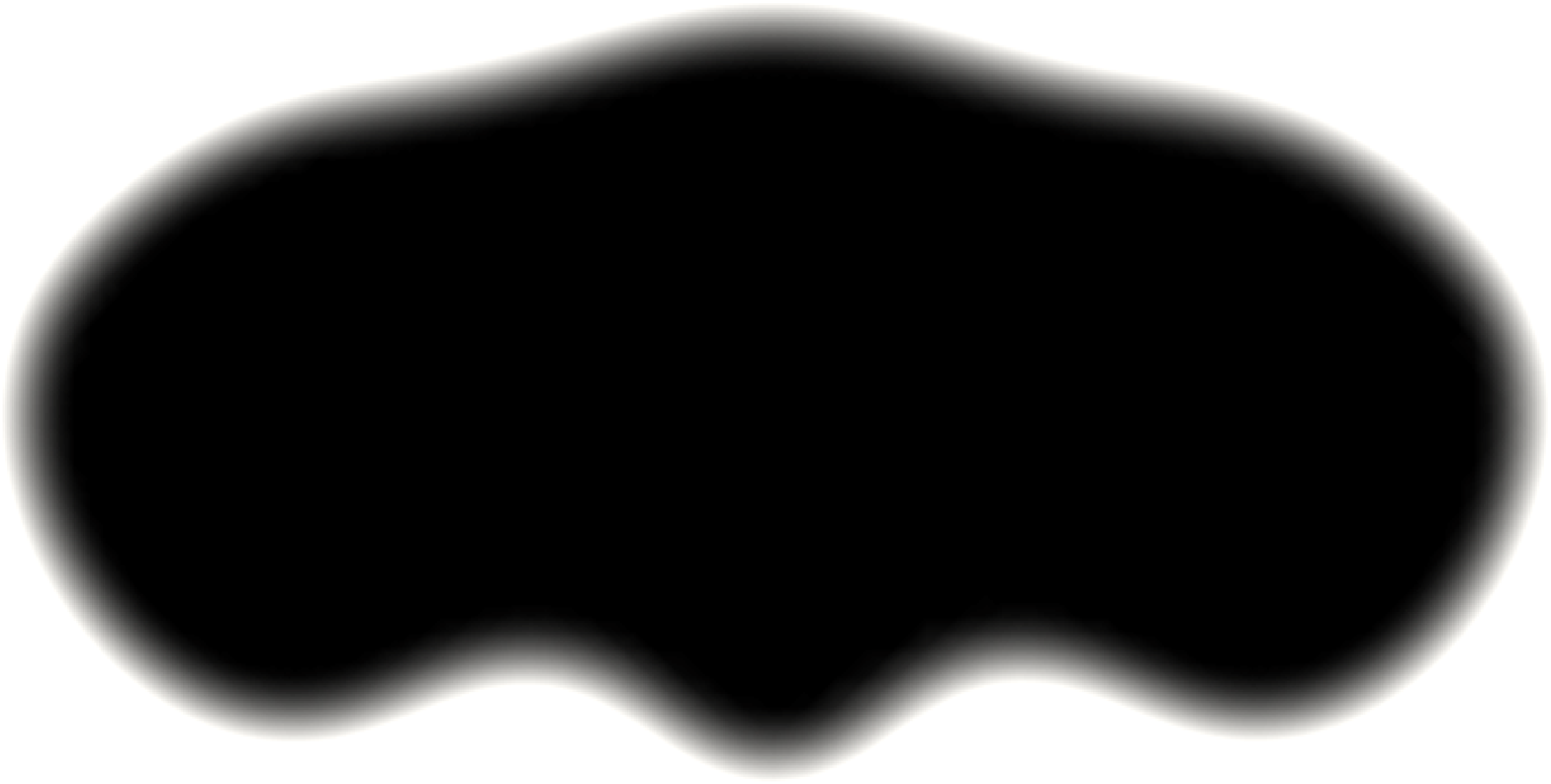}
  }
  \hfill
  \fbox
  {
  \includegraphics[width=0.15\textwidth]{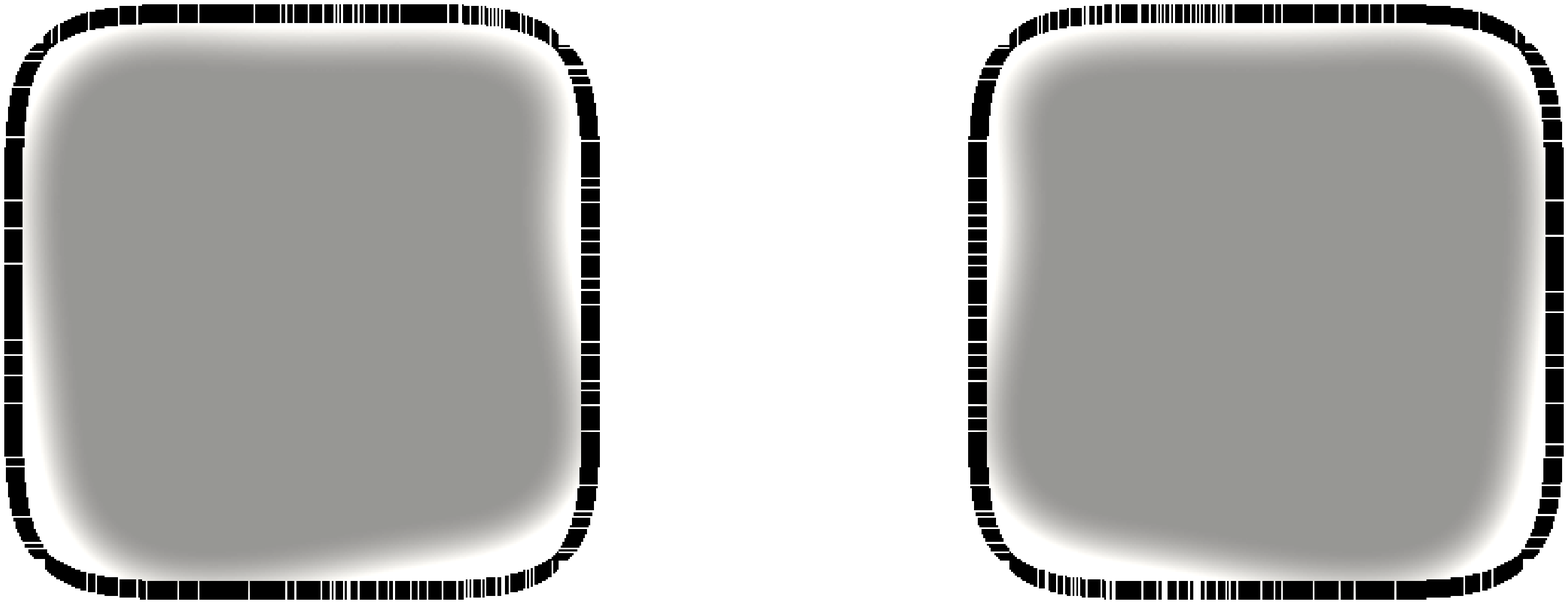}
  }
  \hfill
  
  \caption{The evolution of the phase field $\varphi$.}
  \label{p13fig:num:c2sX:results}
  
\end{figure}

In Figure \ref{p13fig:num:c2sX:results} we depict the temporal evolution
of the phase field $\varphi$ corresponding to the optimal solution at the times $t = 0.00, 0.25,  0.50, 0.75, 1.00$.
The figure additionally includes the zero level line of the desired shape $\varphi_d$ for $t=1.00$.

\begin{figure}
  \centering
  \fbox{
  \includegraphics[width=0.4\textwidth]{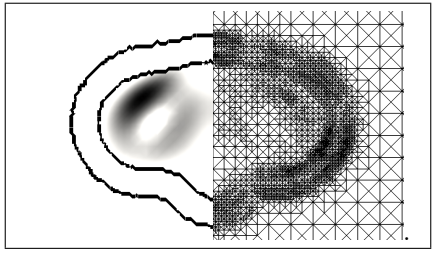}.
  }
  \caption{The magnitude of $v$ in grayscale and the isolines $\varphi \equiv \pm 1$ (left), and the associated triangulation (right).}
  \label{p13fig:num:c2sX:yRefine}
\end{figure}

Regarding the mesh adaptation process,
we observe that the cells are mainly refined in the interfacial region and, in particular, at the border of the diffuse interface.
Such a behavior is typical for the numerical simulation of phase field models.
However, since our error estimator
also contains terms from the Navier-Stokes and the adjoint equation, we further
obtain significant mesh adaptations outside of the interface of the phases, which suggests that these errors should not be neglected, e.g. by a simple interface refinement
technique.
In Figure~\ref{p13fig:num:c2sX:yRefine} we depict the
subdomain $\Omega_u = (0,1) \times (0.5,1.0) \subset \Omega$ at $t=0.7$.
On the left we show $|v|$ in grayscale together with the isolines $\varphi\equiv \pm 1$ in black.
On the right we show the corresponding mesh. Note that the mesh is
symmetric with respect to the central line.

\subsection{Bundle-free implicit programming approach}
\label{p13sec:bundle}

The Algorithm \ref{p13n15:alg:alg:overallAlgorithm} can be further enhanced by exploiting the specific structure of the directional derivative of the control-to-state operator.
Hereby, we apply the descent method directly to the problem \eqref{p13optprob.Ppsi} or \eqref{p13reducopt} (instead of a regularized version) 
and compute a descent direction of $\overline{\mathcal{J}}$ at $u^*$ with $(v^*,\varphi^*,\mu^*)=S(u^*)$ by solving the optimization problem
\begin{align}
\begin{aligned}
\min_{h\in L^2(\Omega;\mathbb{R}^N)^{K-1}} &\overline{\mathcal{J}}'[u^*]( {h})+\| {h}\|^2=(\varphi^*-\varphi_d, {q})+\xi(u^*, {h})+\| {h}\|^2,\\
&\textnormal{s.t. } DS_\Psi[u^*]( {h})=( {q}, {w},{\zeta}),\label{p13direcopt}
\end{aligned}
\end{align}
where the stabilizing term $\| {h}\|^2_{L^2}$ ensures the existence of solutions.
If a solution $h$ of \eqref{p13direcopt} equals zero, then $u^*$ is a B-stationary point, otherwise it is indeed a descent direction, since $\overline{\mathcal{J}}'[u^*]({h})\leq -\|h\|^2< 0$.
In combination with a classical line search procedure, this leads to the following Algorithm \ref{p13descentalgo}.

\begin{algorithm}[H]
  \KwData{Initial data: $\varphi_{a},v_a,u_0$;}
\Repeat{$h_k\leq \epsilon_{tol}$}
{Calculate a descent direction $h_k$ by solving \eqref{p13direcopt}\;
Find a step size $\tau_k$ and a new iterate $u_{k+1}:=u_k+\tau_k h_k$ by performing an Armijo line search along $h_k$ \label{p13line:linesearch}\; 
Set $k:=k+1$.}
 \caption{The descent method for \eqref{p13optprob.Ppsi}}
\label{p13descentalgo}
\end{algorithm}

The convergence of Algorithm \ref{p13descentalgo} is ensured based on the arguments of \cite{Hintermueller2012a}.

\begin{thm}\label{p13convalgo}
 The conceptual Algorithm \ref{p13descentalgo} terminates after finitely many steps for any starting point $u_0$ if either 
  $\tau_k\geq \underline{\tau}>0$ for every $k\in\mathbb{N}$, or $\tau_k\rightarrow 0$ and
  \begin{align}
   \limsup_{k\rightarrow\infty}\frac{\overline{\mathcal{J}}(u_k+\overline{\tau}_k h_k)-\overline{\mathcal{J}}(u_k)-\overline{\tau}_k\overline{\mathcal{J}}'[u_k]({h_k})}{\overline{\tau}_k}\leq 0,
  \end{align}
  where $\overline{\tau}_k>0$ represents the smallest step size for which the line search still fails at step $k$.
\end{thm}

\vspace{2ex}
Motivated by Theorem \ref{p13convalgo}, we include an additional robustification step by performing one step of the penalization algorithm of Subsection \ref{p13sec:penalgo}, if the step size tends to zero.
Thus, the resulting algorithm targets strong stationary points of \eqref{p13optprob.Ppsi}, while guaranteeing at least C-stationarity of the computed solutions.

In order to solve the problem \eqref{p13direcopt}, we take advantage of the fact that it corresponds to a quadratic program, if strict complementarity holds, i.e. if the biactive set associated with the variational inequality \eqref{p13firsttim2} is empty.
Otherwise we employ a regularization of the lower-level problem associated with \eqref{p13direcopt}.

As in the previous subsection, we utilize Taylor-Hood finite elements for the spacial discretization, and supplement the algorithm with a similar adaptive mesh refinement strategy.
Moreover, we solve the discretized CHNS system via a primal-dual active set method.

\begin{figure}
  
  \centering
  
  \includegraphics[width=0.2\textwidth]{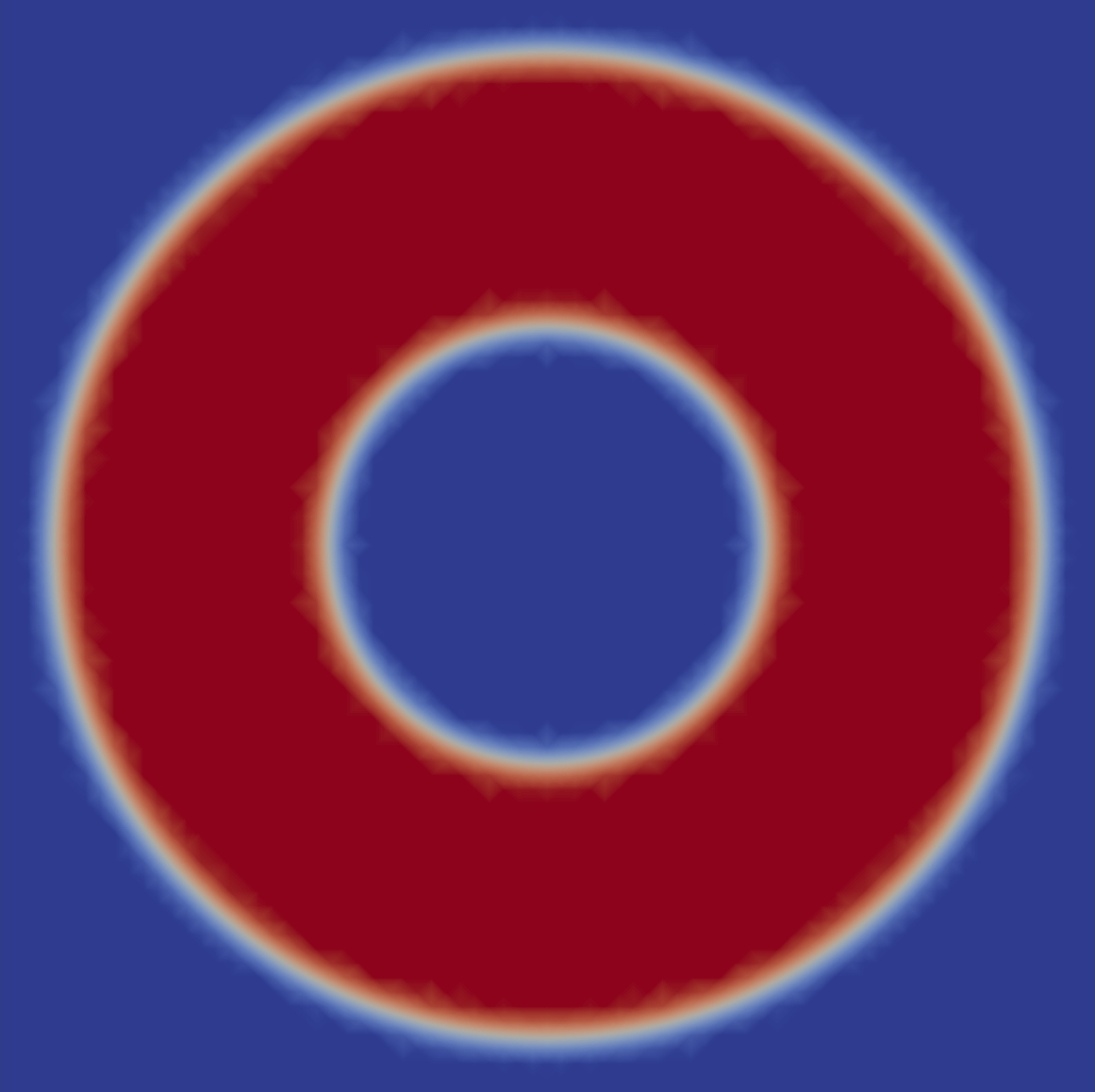}
  \hspace{0.5cm}
  \includegraphics[width=0.2\textwidth]{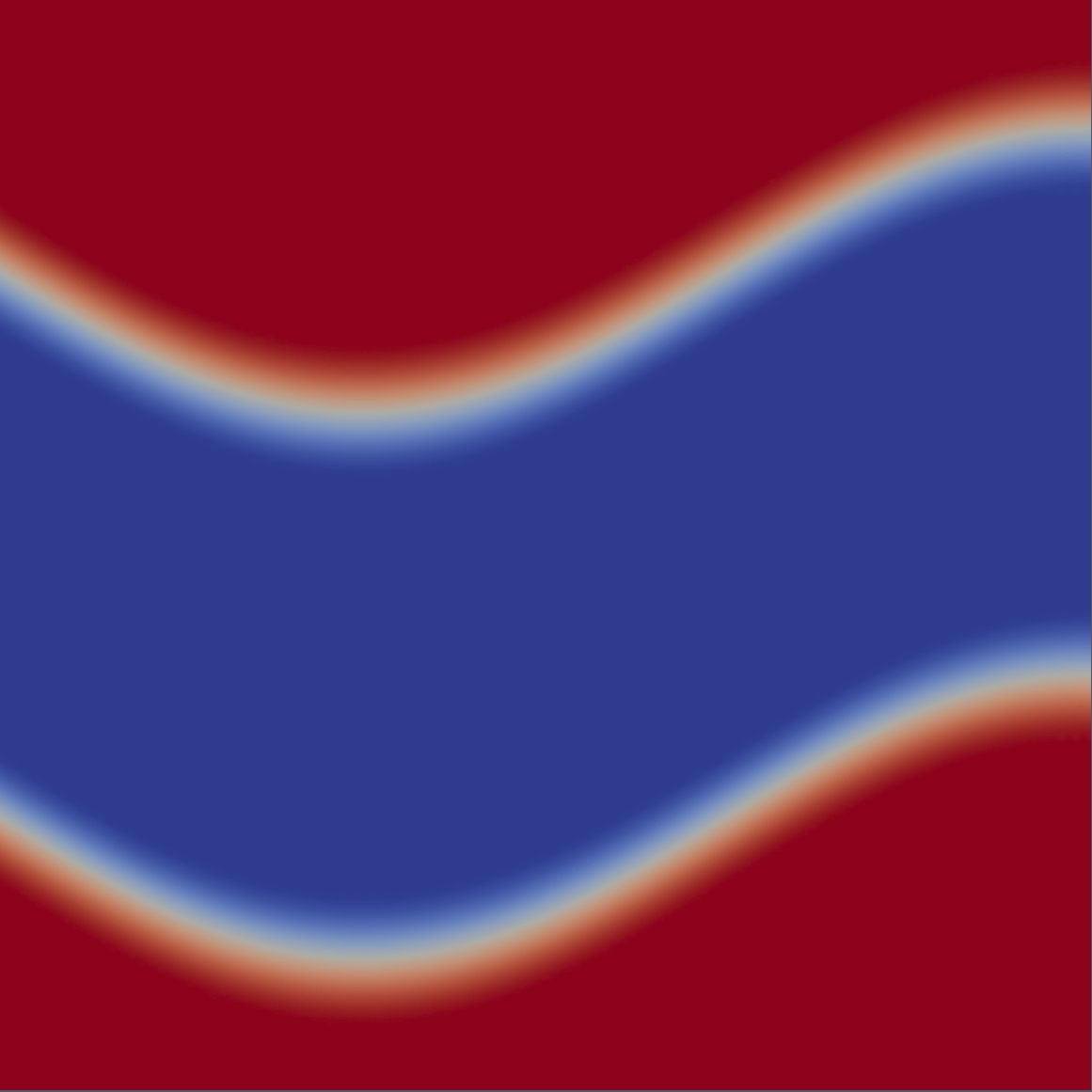}
  \hspace{0.5cm}
  \includegraphics[width=0.2\textwidth]{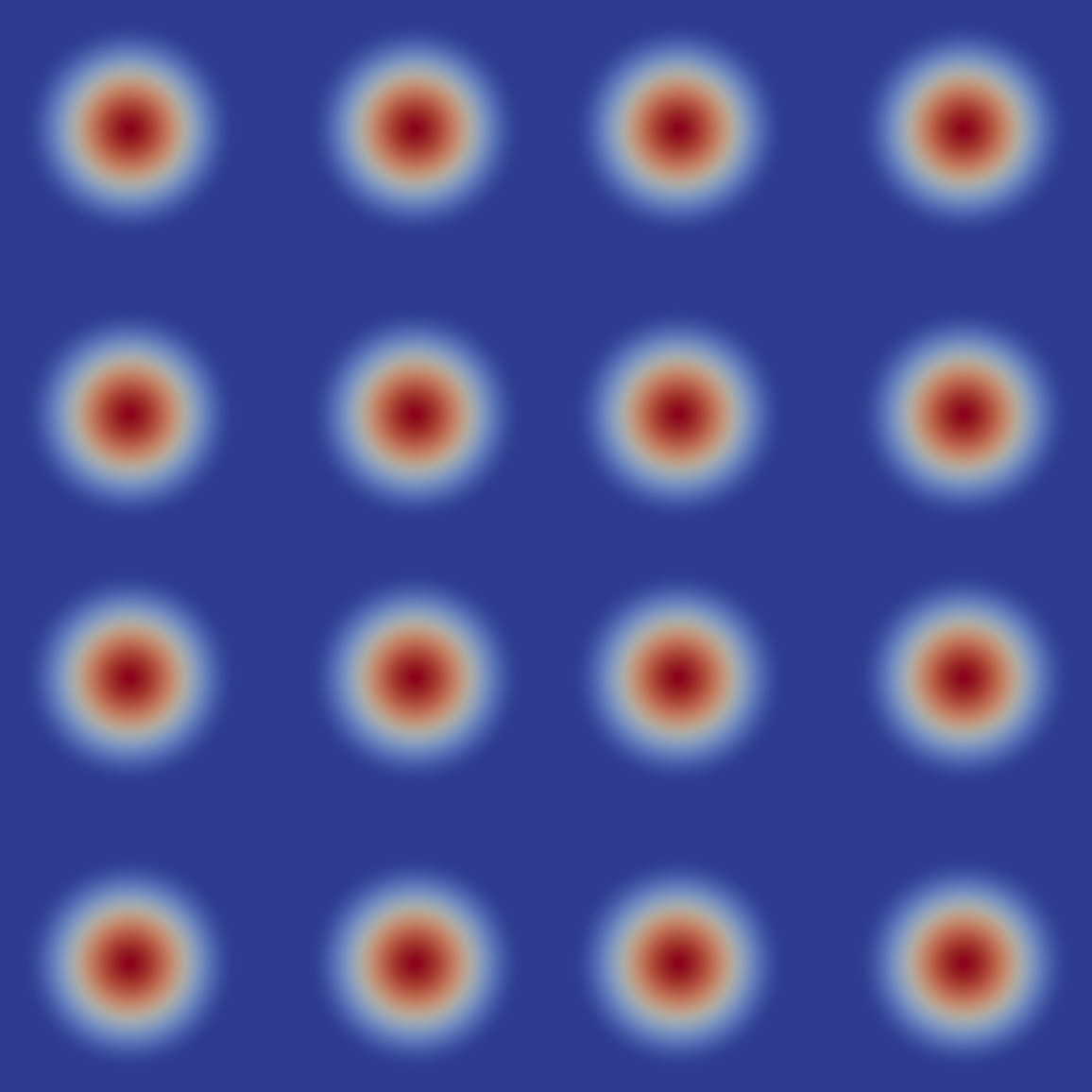}

  \caption{The initial shape $\varphi_0$, the desired shape $\varphi_d$, the ansatz for the control $u$.}
  \label{p13fig:num:R2T}
\end{figure}

In the following example, we aim to transform a ring-shaped initial region into a curved tube, see Figure \ref{p13fig:num:R2T}.
As seen on the right picture, the control acts via 16 locally supported Ansatz functions.

The parameters for the physical model and the adaptation procedure are adopted from the previous example.
In this example, the algorithm terminates at a C-stationary point after performing the Armijo line search (in line \ref{p13line:linesearch}) 276 times.
The maximum number of cells is exceeded after 6 mesh refinement steps.

\begin{figure}
  \centering
  \hfill
      
  \includegraphics[width=0.15\textwidth]{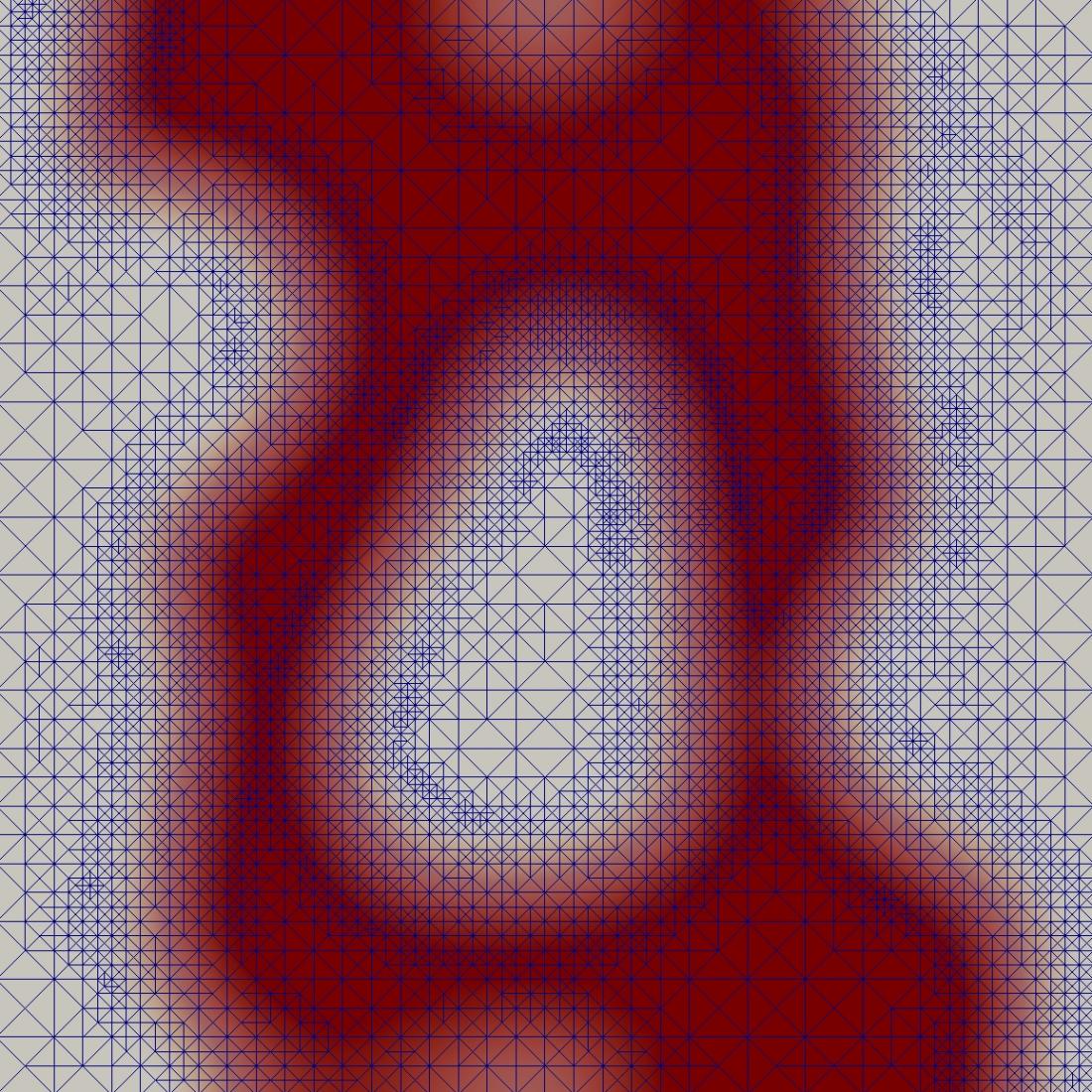}
  \hfill
  \includegraphics[width=0.15\textwidth]{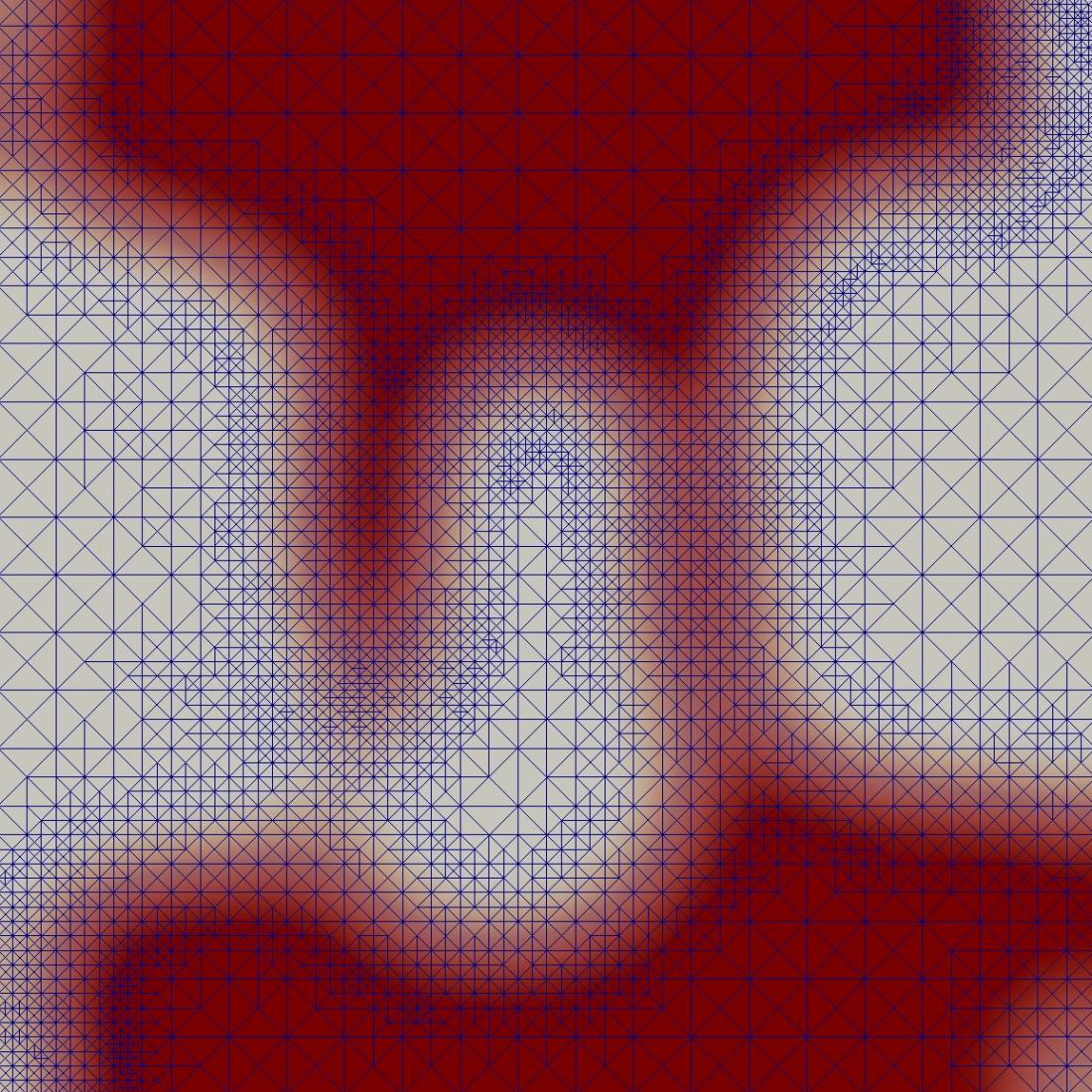}
  \hfill
    \includegraphics[width=0.15\textwidth]{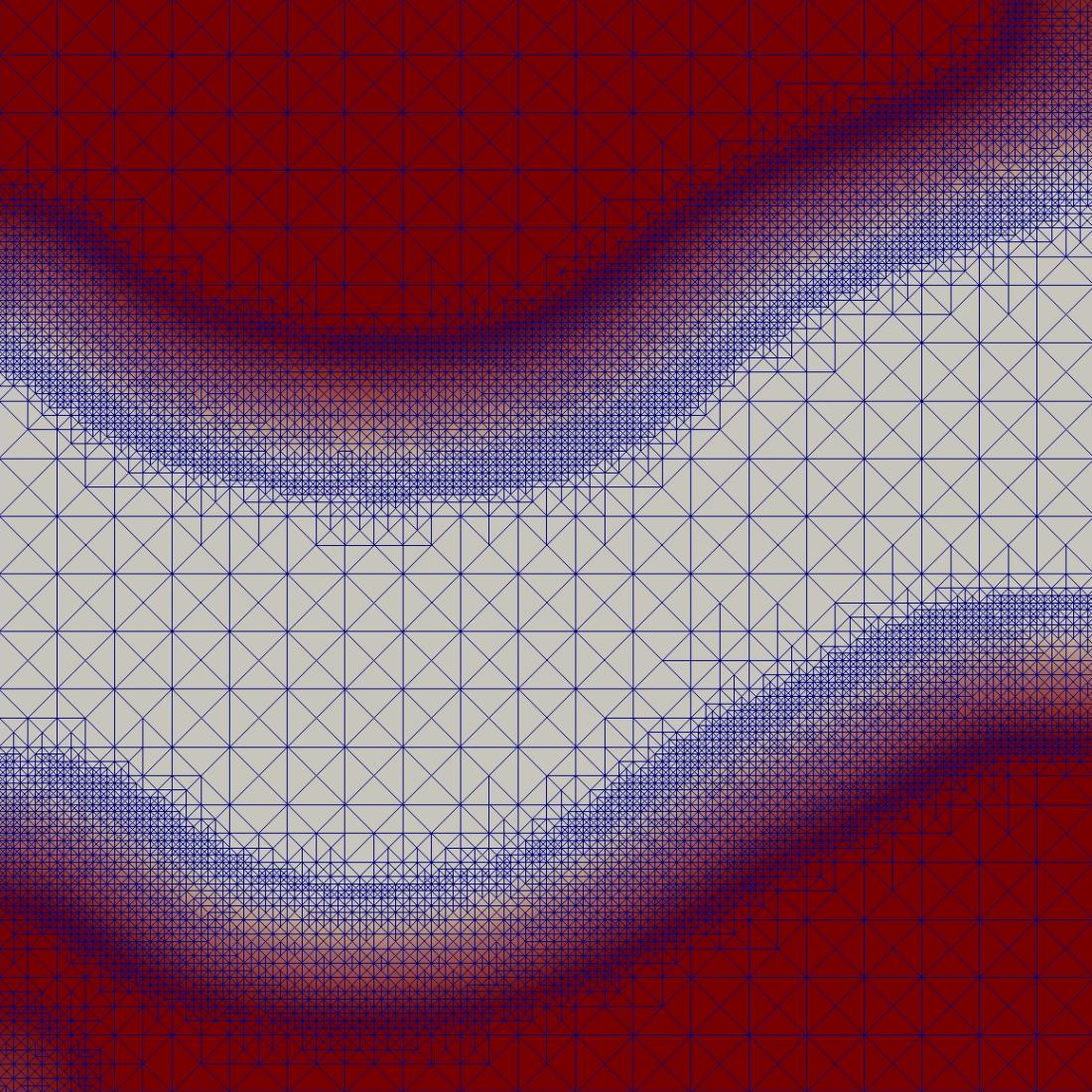}
  \hfill
  \includegraphics[width=0.15\textwidth]{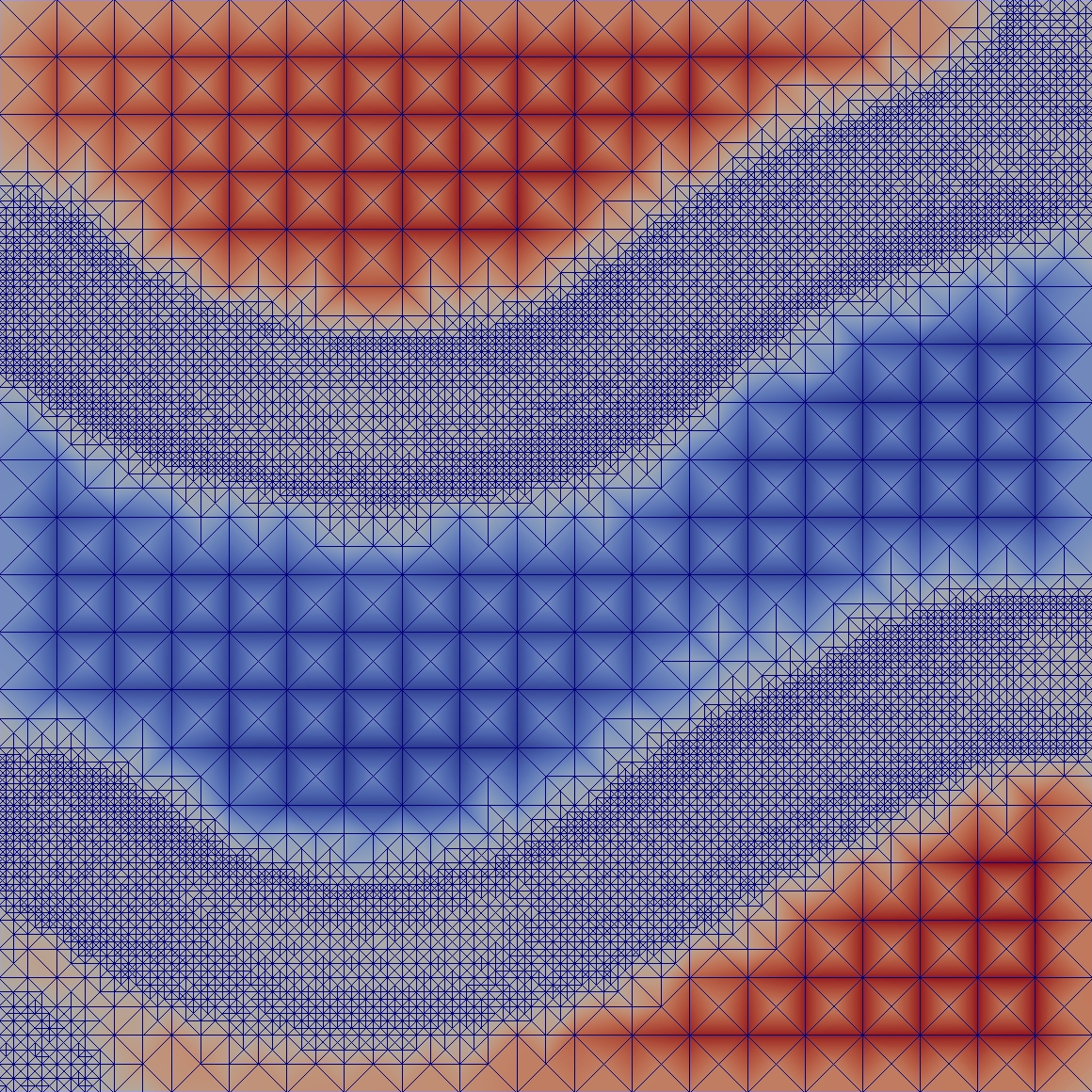}
  \hfill
  \includegraphics[width=0.15\textwidth]{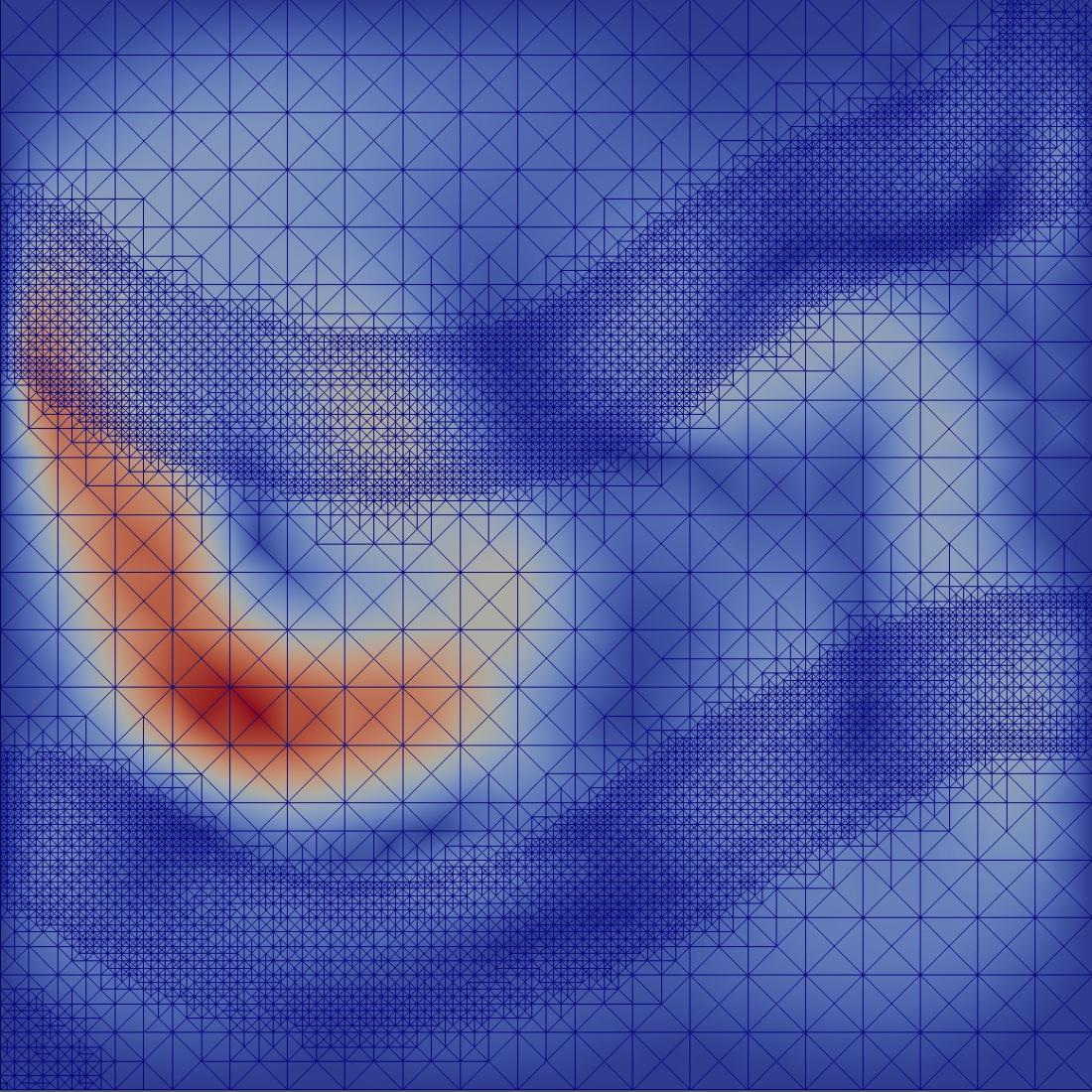}
  \hfill
  
  \caption{The evolution of the phase field $\varphi$, the slack variable $a$ and the magnitude of $v$ at the final time.}
  \label{p13fig:num:R2T:results}
  
\end{figure}

Figure \ref{p13fig:num:R2T:results} presents the computed evolution of the phase field $\varphi$ at the optimal solution along with the associated slack variable $a$ emerging from the primal-dual active set method at the final time.
In addition, we portray the magnitude of the velocity and the underlying mesh at final time.

\section{Model Order Reduction with Proper Orthogonal Decomposition}\label{p13sec:mor}
From a numerical point of view, the simulation and in particular the optimal control of the coupled Cahn-Hilliard Navier-Stokes system \eqref{p13eq:CHNS} are computationally demanding tasks. Although the use of adaptive finite element discretization concepts (see e.g.\ \cite{HHKK}) makes numerical implementation feasible (in comparison to the use of a very fine, uniform discretization), the computational costs can be very large. For this reason, we apply model order reduction using Proper Orthogonal Decomposition (POD-MOR) in order to speed up computation times while ensuring a good approximation quality.\\
In order to construct a low-dimensional surrogate model, the usual POD framework first requires a so-called offline phase, in which high-fidelity solutions (snapshots) of the underlying dynamical system are generated by e.g.\ finite element simulations. From this snapshot set, the POD method finds a proper basis representation of the most relevant information encoded in the snapshots by computing a truncated singular value decomposition or by solving an associated eigenvalue problem. If the snapshots are discretized adaptively in space, the challenge arises that the snapshots are vectors of different lengths due to the different spatial resolutions at each time instance. This does not fit into the standard POD framework which assumes snapshots of the same length.\\
This section is concerned with POD reduced-order modeling using space-adapted snapshots. Section \ref{p13sec:PODHilb} describes the idea to consider the setting from an infinite-dimensional perspective which allows a broad spectrum of discretizations for the snapshots. Then, we derive a POD reduced-order model for the Cahn-Hilliard equations using space-adapted snapshots in Section \ref{p13sec:PODCH} and present numerical results. Moreover, in Section \ref{p13sec:PODNaSt} we consider POD-MOR with space-adapted snapshots for incompressible flow governed by the Navier-Stokes equations, where two strategies are proposed in order to ensure stability of the reduced-order model.

\subsection{POD in Hilbert spaces with space-adapted snapshots}\label{p13sec:PODHilb}
For a comprehensive study of the infinite-dimensional perspective on POD in a Hilbert space setting let us refer to \cite{KV02}, for example. Here, we recall main aspects and provide a practical implementation which is proposed in \cite{GH18}.\\
Let $\{y_h^0, \dots, y_h^{K-1}\} \subset X$ be a given set of snapshots, where $X$ denotes a real, separable Hilbert space and $y_h^i$ for $i=0, \dots, K-1$ are high-fidelity adapted finite element solutions of the underlying dynamical system at different time instances. In particular, each of the snapshots belongs to a different discrete Galerkin space $y_h^i \in V_h^i$ with $V_h^0, \dots, V_h^{K-1} \subset X$. Then, a POD basis of rank $\ell$ is constructed by solving the following equality constrained minimization problem:
\begin{equation}\label{p13eq:POD}
  \min_{\psi_1, \dots, \psi_\ell \in X}  \sum_{i=0}^{K-1} \alpha_j \left\| y_h^i - \sum_{j=1}^\ell (y_h^i,\psi_j)_X \psi_j \right\|_X^2 \quad \text{s.t.} \quad (\psi_i, \psi_j)_X = \delta_{ij} \text{ for } 1 \leq i,j \leq \ell,
\end{equation}
where $\alpha_i$ for $j=0, \dots, K-1$ denote nonnegative weights and $\delta_{ij}$ is the Kronecker symbol. Since the snapshots are spatially adapted, the number of degrees of freedom and/or the location of the node points might differ such that it is not possible to build a corresponding snapshot matrix containing the finite element Galerkin coefficients. For this reason, we assemble the snapshot Gramian defined by
$$\mathcal{K} \in \mathbb{R}^{K\times K}, \quad \mathcal{K}_{ij} := \sqrt{\alpha_i \alpha_j} (y_h^i,y_h^j)_X$$ 
for $i,j=0,\dots,K-1$. In order to set up the matrix $\mathcal{K}$, we only require that the snapshots belong to the same Hilbert space $X$ in order to evaluate the inner product $(\cdot,\cdot)_X$. Solving an eigenvalue problem for $\mathcal{K}$, i.e.\
$$ \mathcal{K}\upphi_i = \lambda_i \upphi_i \quad \text{for } i = 1, \dots, \ell$$
delivers eigenvalues $\lambda_1 \geq \dots \geq \lambda_\ell \geq 0$ and eigenvectors $\{\upphi_1, \dots, \upphi_\ell\} \subset \mathbb{R}^{K}$ which suffice to set up the POD reduced-order model, see \cite[Section 4]{GH18} for more details. The advantage of this perspective is that it allows a broad spectrum of discretization techniques and includes the case of $r$-adaptivity, for example. However, in this case the evaluation of the inner products $(y_h^i,y_h^j)_X$ might get involved such that the necessity of e.g.\ parallelization becomes evident for practical implementations. In case of $h$-adapted snapshots using hierarchical, nested meshes, it is reasonable to express the snapshots with respect to a common finite element space as proposed in \cite{URL16}.

\subsection{POD reduced-order modeling for the Cahn-Hilliard system}\label{p13sec:PODCH} 

Let us consider the weak formulation of the Cahn-Hilliard equations \eqref{p13eq:CH1}-\eqref{p13eq:CH2} with boundary conditions \eqref{p13eq:bdryCH} and an initial condition for the phase field \eqref{p13eq:initialCH}, where we assume the velocity $v$ to be given and fixed. The weak form reads as: Find a phase field $\varphi \in W(0,T;H^1(\Omega))$ with $\varphi|_{t=0} = \varphi_a$ and a chemical potential $\mu \in L^2(0,T;H^1(\Omega))$ such that
 \begin{subequations}\label{p13eq:CHweak}
 \begin{alignat}{4}
 \frac{d}{dt}(\varphi(t),\phi)_{L^2(\Omega)} +(v\nabla\varphi(t),\phi)_{L^2(\Omega)} + m(\nabla \mu(t), \nabla \phi)_{L^2(\Omega)} & = && \; 0 && \quad \forall \phi \in H^1(\Omega), \label{p13eq:CH1weak} \\
 \sigma \epsilon (\nabla \varphi(t), \nabla \phi)_{L^2(\Omega)} +\frac{\sigma}{\epsilon}(\Psi_0'(\varphi(t)) -\kappa\varphi(t),\phi)_{L^2(\Omega)} - (\mu(t),\phi)_{L^2(\Omega)} & = && \; 0 && \quad \forall \phi \in H^1(\Omega). \label{p13eq:CH2weak} 
\end{alignat}
\end{subequations}
Note that in \eqref{p13eq:CHweak} we assume for simplicity a constant mobility $m>0$ and sufficient regularity for $\Psi_0$. In order to derive an associated POD reduced-order model, we approximate the phase field $\varphi$ and the chemical potential $\mu$ by a POD Galerkin ansatz given as $\varphi(t) \approx \varphi_\ell(t) = \sum_{j=1}^\ell c_j(t) \psi_j$ and $\mu(t) \approx \mu_\ell(t) = \sum_{j=1}^\ell w_j(t) \psi_j$. In \cite{GH18,GHpamm}, we construct separate POD reduced spaces for the phase field and the chemical potential, respectively. In contrast, here we compute the POD modes $\psi_j$ for $j=1,\dots,\ell$ according to \eqref{p13eq:POD} from space-adapted finite element snapshots of the phase field and use the same POD modes in the Galerkin ansatz for both phase field and chemical potential. Using the POD space $V_\ell = \textnormal{span}\{\psi_1, \dots, \psi_\ell\} \subset H^1(\Omega)$ as trial and test space leads to the following POD reduced-order model for the Cahn-Hilliard equations: Find a phase field $\varphi_\ell \in V_\ell$ with $\varphi|_{t=0} = \mathcal{P}_\ell \varphi_a$ and a chemical potential $\mu_\ell \in V_\ell$ such that
 \begin{subequations}\label{p13eq:CHPOD}
 \begin{alignat}{4}
 \frac{d}{dt}(\varphi_\ell(t),\psi)_{L^2(\Omega)} +(v\nabla\varphi_\ell(t),\psi)_{L^2(\Omega)} + m(\nabla \mu_\ell(t), \nabla \psi)_{L^2(\Omega)} & = && \; 0 && \quad \forall \psi \in V_\ell, \label{p13eq:CH1POD} \\
 \sigma \epsilon (\nabla \varphi_\ell(t), \nabla \psi)_{L^2(\Omega)} +\frac{\sigma}{\epsilon}(\Psi_0'(\varphi_\ell(t)) -\kappa\varphi_\ell(t),\psi)_{L^2(\Omega)} - (\mu_\ell(t),\psi)_{L^2(\Omega)} & = && \; 0 && \quad \forall \psi \in V_\ell. \label{p13eq:CH2POD} 
\end{alignat}
\end{subequations}
By $\mathcal{P}_\ell: V \to V_\ell$ we denote the orthogonal projection onto the POD space. Note that in \eqref{p13eq:CHPOD}, the evaluation of the nonlinear term $\Psi_0'(\varphi_\ell(t))$ is dependent on the full-order dimension. The treatment of nonlinearities is a well-known challenge within POD-MOR. In order to enable an efficient evaluation of the nonlinearity which is related to the low-order dimension $\ell$ of the reduced system, a linearization can be considered, compare \cite{GH18} for more details. Alternatively, so-called hyper-reduction methods like EIM \cite{BMNP04}, DEIM \cite{CS10} or DMD \cite{AK17} can be applied.\\

\subsection{Numerical example of POD-MOR for the Cahn-Hilliard system}\label{p13sec:numericspodch}

In this Section, we numerically investigate two major issues within POD-MOR for the Cahn-Hilliard equations: 

\begin{itemize}
 \item[(i)] How does the regularity of the free energy $\Psi_0$ effect the accuracy of the POD reduced-order model?
 \item[(ii)] How does the use of spatial adaptivity in the offline phase for snapshot generation influence the computational times and the accuracy of the POD reduced-order model?
\end{itemize}

\begin{figure}[htbp]
 \includegraphics[scale=0.25]{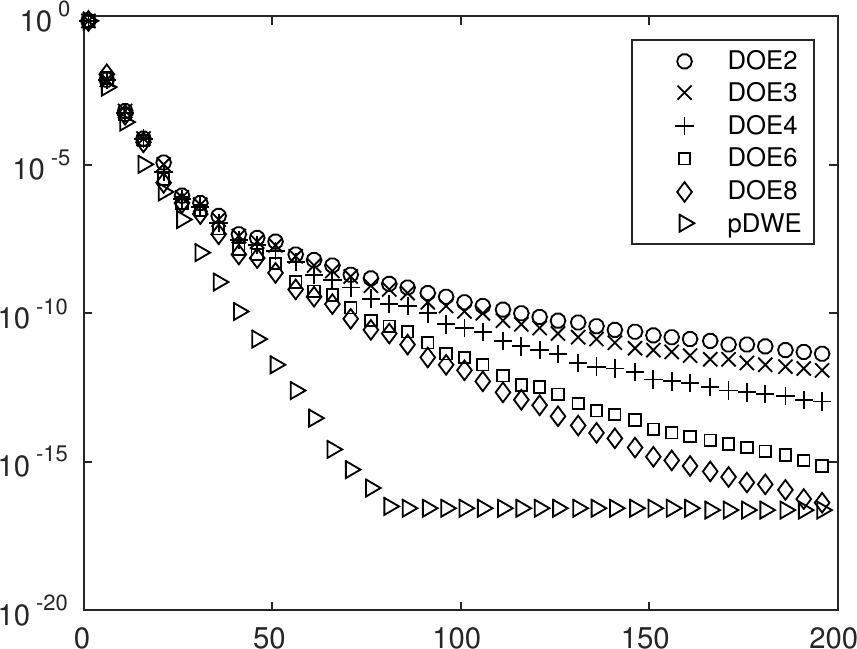}
 \caption{Decay of the normalized eigenvalues for the phase field $\varphi$ considering a Moreau-Yosida relaxation (DOE$r$) for different relaxation parameters $r$ and a polynomial free energy (pDWE)}\label{p13fig:ev}
\end{figure}

\noindent The first aspect (i) is studied numerically in \cite{Alf}. The initial phase field is given as a circle in a two-dimensional domain, which is transported in horizontal direction over time. In this simulation, a uniform and static discretization in space is used to generate the snapshots and a POD basis is computed with respect to the $X=L^2(\Omega)$-inner product. The decay of the normalized eigenvalues is shown in Figure \ref{p13fig:ev}. It compares the use of a smooth double well potential $\Psi_0(\varphi) = \frac{1}{4} \varphi^4$ (pDWE) to the use of a Moreau-Yosida relaxation of the double-obstacle potential given as $\Psi_0(\varphi) = \frac{s}{r}(|\max(0,\varphi-1)|^r + |\min(0,\varphi+1)|^r)$ for different values of $r$ (DOE$r$). We observe that the smoother the considered free energy is, the faster is the decay of the eigenvalues. This is similar to a well-known behavior in Fourier analysis, where the decay of the Fourier coefficients depends on the smoothness of the object. For POD reduced-order modeling this means that if a potential with lower regularity is used, then more POD modes are needed for an adequate approximation than using a smooth potential.

In future research, we plan to apply POD model order reduction for the Cahn-Hilliard equations using a nonsmooth double-obstacle potential. This involves reduced-order modeling for variational inequalities, see e.g.\ \cite{BHSW} for a reduced-order technique for Black-Scholes and Heston models.\\

\begin{figure}[htbp]
 \includegraphics[scale=0.16]{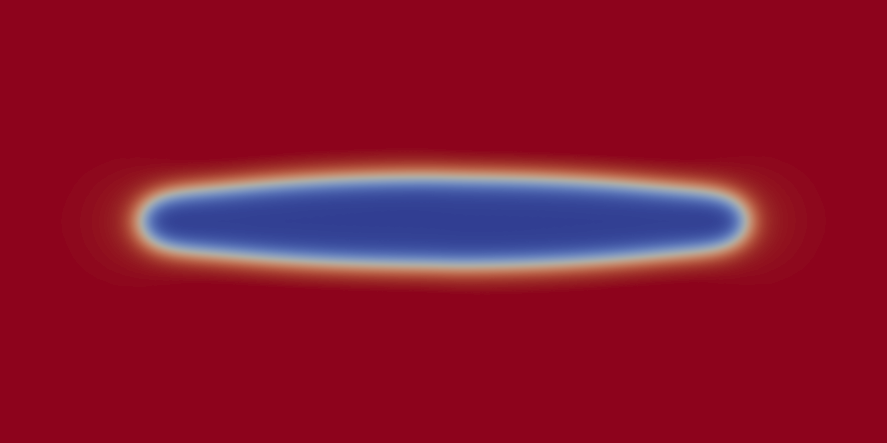}
 \includegraphics[scale=0.16]{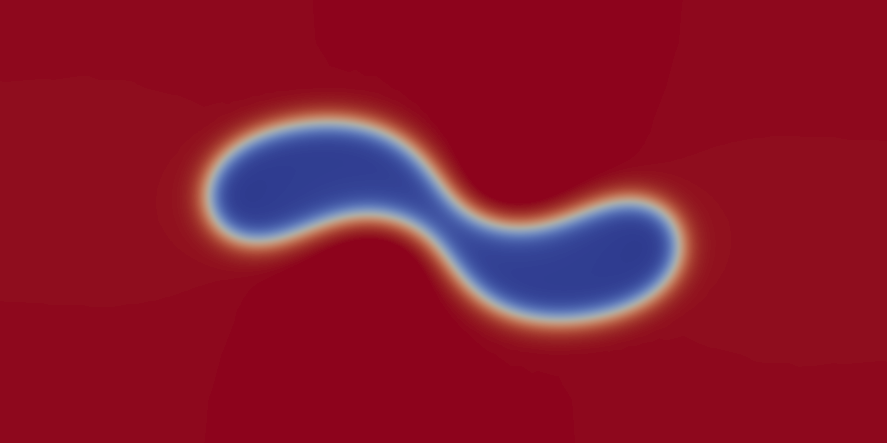}
 \includegraphics[scale=0.16]{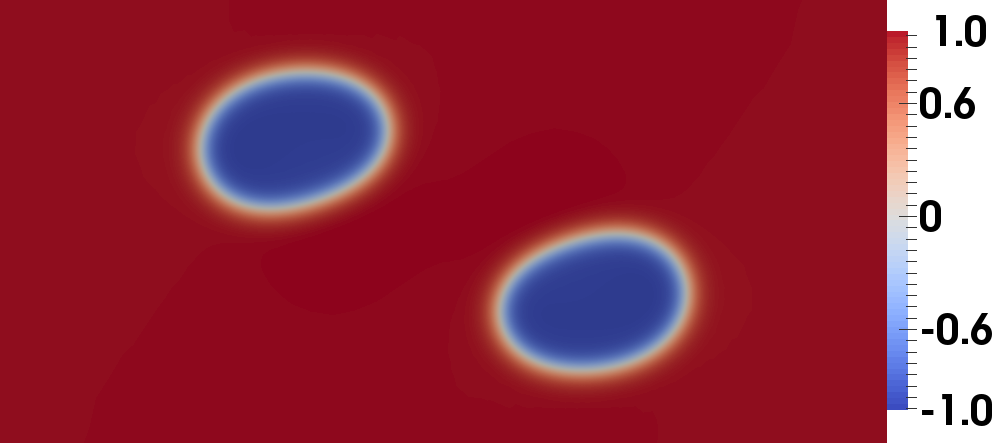}\\[0.1cm]
  \hspace*{-0.6cm} \includegraphics[scale=0.16]{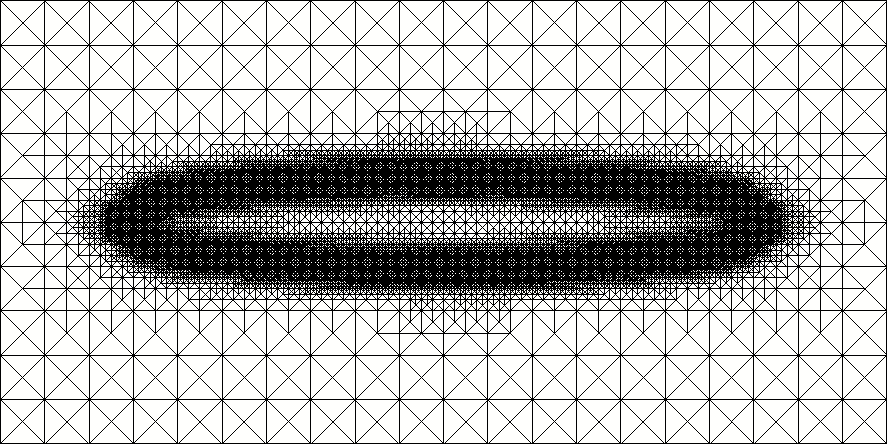} 
 \includegraphics[scale=0.16]{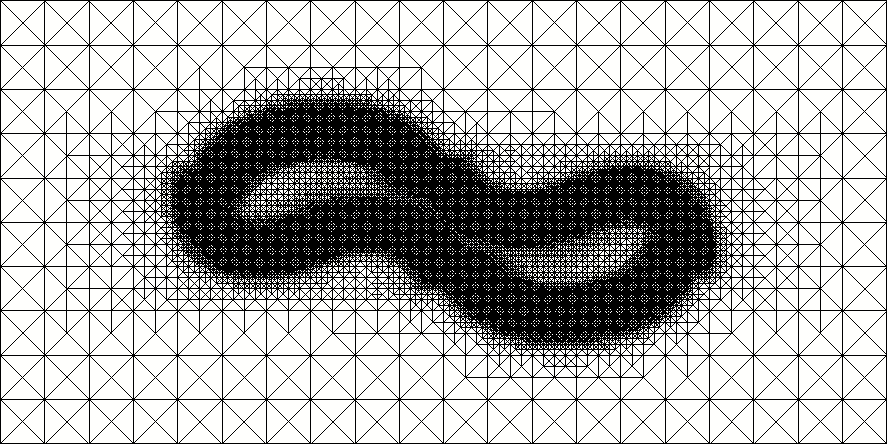}
 \includegraphics[scale=0.16]{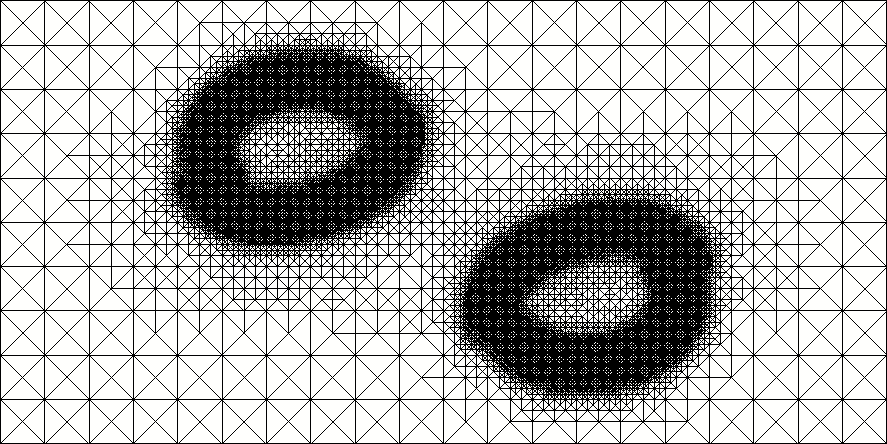}
\caption{Finite element snapshots of the phase field at $t=0$, $t=T/2$ and $t=T$ (top) with the associated adapted finite element meshes (bottom)}\label{p13fig:adfe}
\end{figure}

\noindent For the second aspect (ii), let us consider the following setting: the spatial domain is $\Omega = (0,2)\times(0,1) \subset \mathbb{R}^2$, the mobility is $m=1.0$, the interface parameter is $\epsilon = 0.02$ and the potential $\Psi_0$ is the smooth double-well energy. The initial condition has the shape of an ellipse. We consider a solenoidal velocity field $y=(y_1,y_2)$ given by
$$ y_1(x) = c \sin(\pi x_0) \cdot \cos (\pi x_1), \; y_2(x) = -c \sin(\pi x_1) \cdot \cos(\pi x_0) \text{ for } x_0 \leq 1$$
and
$$y_1(x) = -c \sin(\pi x_0) \cdot \cos (\pi x_1), \; y_2(x) = c \sin(\pi x_1) \cdot \cos(\pi x_0) \text{ for } x_0 > 1,$$
where $x=(x_0,x_1)$. In this example, we choose $c=70$, such that the velocity field leads to a break-up of the ellipse into two separate droplets. This topology change can be handled naturally due to the consideration of a diffuse interface approach.\\
For the temporal discretization, we use an unconditional gradient stable scheme based on a convex-concave splitting of the potential according to \cite{ES93,Eyr98}. As time step size we use $\tau = 2.5\cdot 10^{-5}$ and perform $K=300$ time steps. For the spatial discretization, we use $h$-adapted piecewise linear and continuous finite elements. The solutions to the adaptive finite element simulation at initial time, half time and end time with the associated adapted meshes are shown in Figure \ref{p13fig:adfe}. The number of node points varies between 16779 and 19808 and the finite element simulation time is 1674 sec.

\begin{figure}[htbp]
 \includegraphics[scale=0.16]{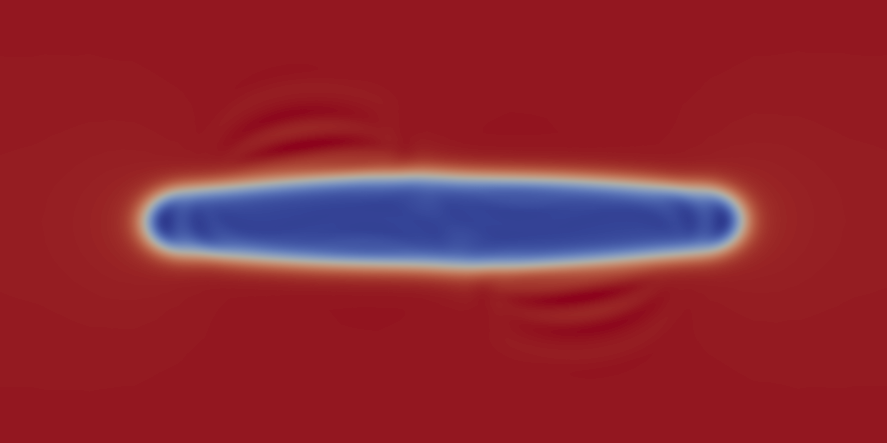}
 \includegraphics[scale=0.16]{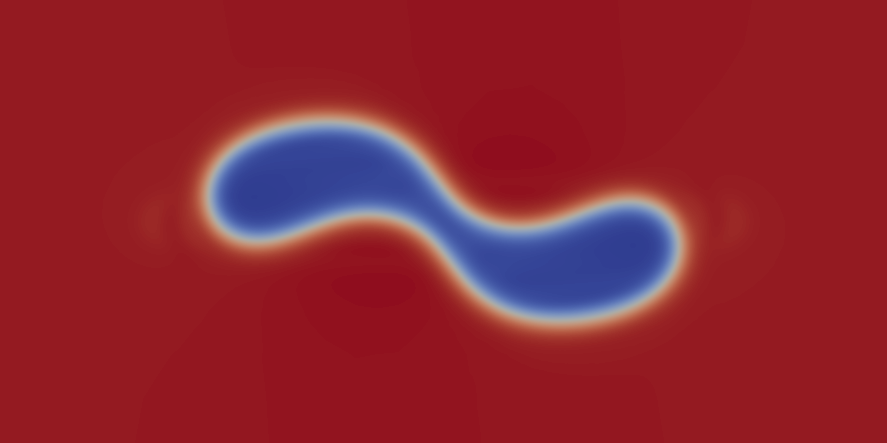}
 \includegraphics[scale=0.16]{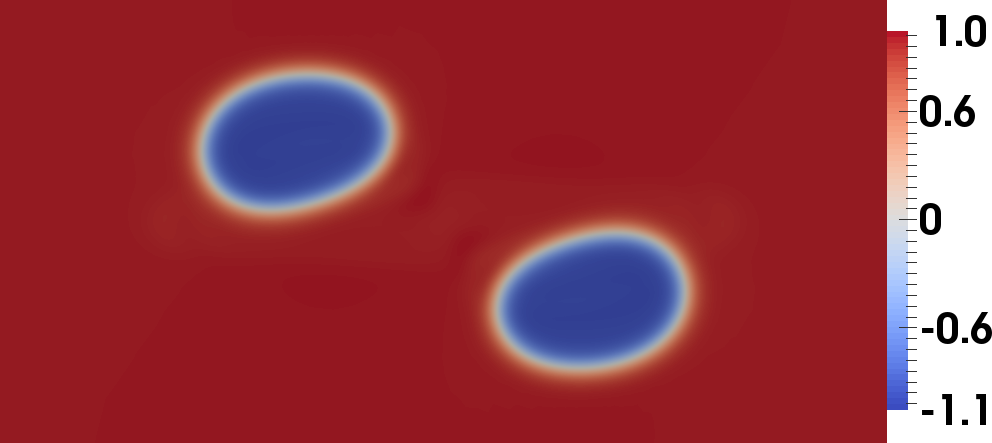}\\[0.1cm]
  \includegraphics[scale=0.16]{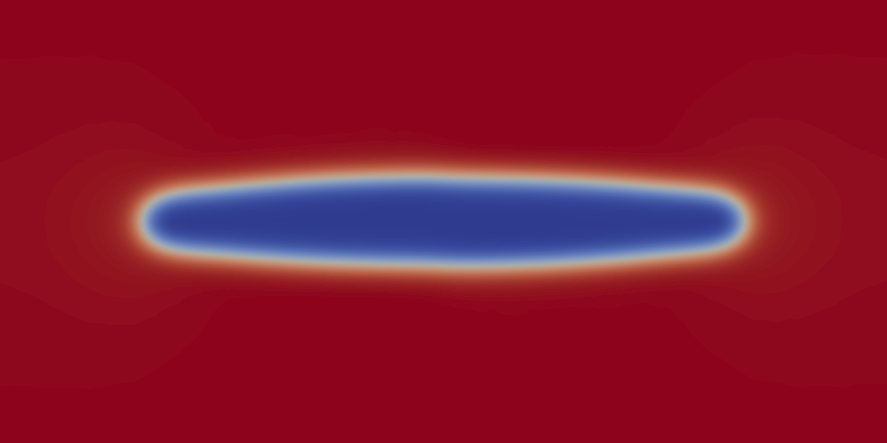}
 \includegraphics[scale=0.16]{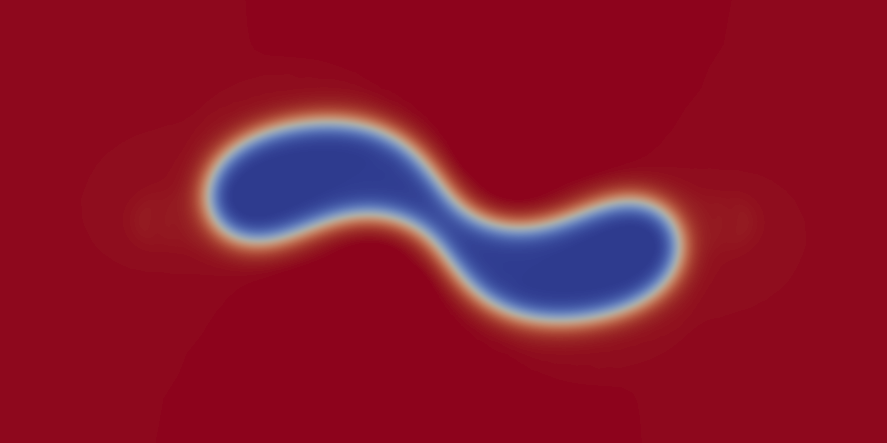}
 \includegraphics[scale=0.16]{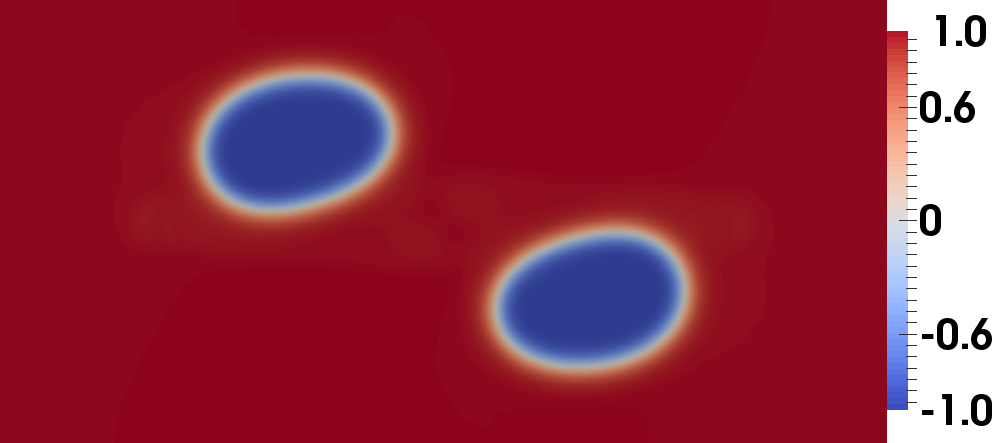}
\caption{POD reduced-order approximation of the phase field at $t=0$, $t=T/2$ and $t=T$ using $\ell=10$ POD modes (top) and $\ell=20$ POD modes (bottom)}\label{p13fig:phirom}
\end{figure}

In order to construct a POD reduced-order model, we utilize the adapted finite element solutions for the phase field as snapshots in \eqref{p13eq:POD}, where we choose $X=L^2(\Omega)$ for the norm and inner products. The resulting solutions for a POD reduced-order model of dimension $\ell=10$ and $\ell=20$ are shown in Figure \ref{p13fig:phirom} at the initial, half and end time. In the approximations using $\ell=10$ POD modes, we observe oscillations due to the transport term, which are smoothened out by enlarging the reduced dimension. We note that POD model order reduction for systems involving a dominant transport is challenging and refer to \cite{CMS,RSSM18,SPO16,WWXI} for different solution concepts.

\noindent The relative $L^2(0,T;\Omega)$-error between the adaptive finite element solution and the POD reduced-order solution using $\ell=20$ POD modes is $2.793 \cdot 10^{-4}$. The solution time for the reduced-order simulation is 88 sec, which leads to a speed up factor of $19$ compared to the time needed for the adaptive finite element simulation. Note that the reduced-order model still depends on the finite element dimension, since an expansion of the reduced solution to the full-order model is needed for the evaluation of the nonlinearity. In order to enable an efficient evaluation of the nonlinearity which is related to the reduced-order dimension, the use of hyper-reduction methods like DEIM is needed. This leads to a further speedup, such that the solution of the reduced-order system takes only a fraction of seconds (compare e.g.\ \cite[Table 5]{GH18}). However, especially in the case of lower regularity of the potential, we observe instabilities. In future research, we plan to derive a stable POD reduced-order model including hyper-reduction for systems with nonlinearities of low regularity. Moreover, we refer to \cite{UK17} for an energy stable model order reduction for the Allen-Cahn equation.\\
For further details on POD with space-adapted snapshots and additional numerical test runs, we refer to \cite{GH18,GHpamm}.\\
The speedup in the computational times when replacing the high-fidelity finite element model by the POD reduced-order surrogate especially pays off in multi-query scenarios like optimal control. In this case, a repeated solution of the associated state and adjoint equations is necessary in order to find a minimum to a given cost functional. We refer to \cite{GHS19} for an optimal control of a Cahn-Hilliard system, where the control enters the equations as velocity in the transport term. A reduced-order model using space-adapted snapshot data is used. A different optimal control problem for the Cahn-Hilliard system is considered in \cite{Alf18}, where the control enters as a right-hand side in \eqref{p13eq:CH1weak}. Within a POD trust-region framework according to \cite{AFS}, the reduced-order model accuracy is evaluated by the Carter condition. This guarantees a relative gradient accuracy and indicates whether an enlargement of the reduced dimension or a POD basis update with space-adapted snapshots at the current optimization iterate is necessary.

\subsection{Stable POD reduced-order modeling for Navier-Stokes with space-adapted snapshots}\label{p13sec:PODNaSt}
Let us now consider the Navier-Stokes system \eqref{p13eq:NaSt1}-\eqref{p13eq:NaSt2} for a single-phase system in strong form, i.e.\ 
\begin{subequations}\label{p13eq:NaStstrong}
\begin{alignat}{4}
  \partial_t v + (v \cdot \nabla) v - \frac{1}{Re} \Delta v + \nabla p & = & f && \quad \text{in } (0,T)\times\Omega, \\
 \textnormal{div} \; v & = & 0 && \quad \text{in } (0,T)\times\Omega,\label{p13eq:NaStstrong2}
\end{alignat}
\end{subequations}
equipped with homogeneous Dirichlet boundary conditions $v = 0$ on $\partial \Omega$ and an initial condition for the velocity \eqref{p13eq:initialNaSt}. In order to derive a fully discrete formulation of \eqref{p13eq:NaStstrong}, we first discretize in time using an implicit Euler scheme, which allows to use a different (adaptive) finite element space at each time instance. Let $t_0 = 0 < t_1 < \dots t_{K-1}=T$ denote a time grid with constant time step size $\tau$ and let $(V_h^i,Q_h^i)$ for $i=0,\dots, K-1$ denote inf-sup stable Taylor-Hood finite element pairs. Then, the fully discrete Navier-Stokes systems reads as: for given $v_h^0=v_a$ find $v_h^1 \in V_h^1, \dots, v_h^{K-1} \in V_h^{K-1}$ and $p_h^1 \in Q_h^1, \dots, p_h^{K-1} \in Q_h^{K-1}$ such that 
\begin{subequations}\label{p13eq:NaStweak}
\begin{alignat}{4}
  \left( \frac{v_h^{i}-v_h^{i-1}}{\tau},w \right) + ((v_h^i \cdot \nabla) v_h^i,w) + \frac{1}{Re} (\nabla v_h^i,\nabla w) + b(w,p_h^i) & = && \; \langle f(t_i),w \rangle && \quad \forall w \in V_h^i, \\
 b(v_h^i,q) & = && \; 0 && \quad \forall q \in Q_h^i,\label{p13eq:NaStweak2}
\end{alignat}
\end{subequations}
for $i=1,\dots,K-1$, where $(\cdot, \cdot)$ denotes the $L^2(\Omega)$-inner product and $\langle \cdot, \cdot \rangle$ is the duality pairing of $H_0^1(\Omega)$ with $H^{-1}(\Omega)$. Moreover, we introduce $b(w,q) := -(q,\nabla \cdot v)$ such that the strong divergence-free condition \eqref{p13eq:NaStstrong2} is now postulated in a weak form in \eqref{p13eq:NaStweak2}. In order to derive the POD reduced-order model, we compute a POD basis from the space-adapted solutions from \eqref{p13eq:NaStweak} according to Section \ref{p13sec:PODHilb}. In particular, we introduce reduced spaces $V_\ell$ and $Q_\ell$ for the velocity and pressure and search for reduced approximations $\{v_\ell^1,\dots,v_\ell^{K-1}\} \in V_\ell$ and $\{p_\ell^1,\dots,p_\ell^{K-1}\}$ such that
\begin{subequations}\label{p13eq:NaStpod}
\begin{alignat}{4}
  \left( \frac{v_\ell^{i}-v_\ell^{i-1}}{\Delta t},w \right) + ((v_\ell^i \cdot \nabla) v_\ell^i,w) + \frac{1}{Re} (\nabla v_\ell^i,\nabla w) + b(w,p_\ell^i) & = && \; \langle f(t_i),w \rangle && \quad \forall w \in V_\ell, \\
 b(v_\ell^i,q) & = && \; 0 && \quad \forall q \in Q_\ell.
\end{alignat}
\end{subequations}
The difficulty consists in the fact that stability of \eqref{p13eq:NaStpod} is not ensured for all choices of $(V_\ell,Q_\ell)$. For this reason, in \cite{GHLU18} we provide two solution concepts:
\begin{itemize}
 \item[(i)] A velocity ROM in the spirit of \cite{Sir} using an optimal projection onto a weak divergence-free space,
 \item[(ii)] A velocity-pressure ROM using a supremizer stabilization technique in the spirit of \cite{BMQR,RV}.
\end{itemize}
In the first approach (i), we utilize the following optimal projection. For a given function $v \in X$ find a reference function $\tilde{v}$ in a reference velocity function space $\tilde{V}$ such that it fulfills 
$$  \min_{u\in\tilde{V}}\frac{1}{2}\|v-u\|_{X}^2\quad\text{s.t.}\quad b(u,q)=0\quad\;\forall q\in\tilde{Q}.$$
This projection is computed either for each of the space-adapted velocity snapshots $\{v_h^1,\dots,v_h^{K-1}\}$ or for each of the velocity POD basis functions $\{\psi_1^v,\dots,\psi_\ell^v\}$ computed from velocity snapshots according to \eqref{p13eq:POD}. Then, a common weak divergence-free property is inherited in the reduced-order model, which leads to a cancellation of the pressure term and continuity equation from \eqref{p13eq:NaStpod}, such that the reduced system is stable by construction. Particular attention must be paid to the treatment of inhomogeneous boundary conditions, for which we refer to \cite[Section 6]{GHLU18} for details.\\
The second approach (ii) utilizes a supremizer enrichment technique. After computing separate POD bases $\{\psi_1^v,\dots,\psi_\ell^v\}$ and $\{\psi_1^p,\dots,\psi_\ell^p\}$ for the velocity and pressure, respectively, we enrich the reduced velocity space by stabilization functions. These are computed as follows: for a given $q \in L_0^2(\Omega)$ find $\mathbb{T}q \in \tilde{V}$ such that
$$ (\mathbb{T}q,\phi)_{H_0^1(\Omega)} = b(\phi,q) \quad \forall \phi \in \tilde{V}.$$
Then, as supremizer functions, we choose $\{\mathbb{T}\psi_1^p,\dots,\mathbb{T}\psi_\ell^p\}$. The inf-sup stability of the resulting velocity-pressure reduced-order model follows from the inf-sup stability of the finite element model, see \cite[Section 5.2]{GHLU18} for the proof.

\section{Outlook}
\label{p13sec:out}
   
In the second phase of the priority programme 1962 we consider shape optimization with instationary fluid flow in a diffuse interface setting. We will provide a well-posed formulation for shape optimization in instationary fluids with general cost functionals, which on the one hand allows for topological changes and imposes no geometric constraints on the optimal shape, and on the other hand overcomes some potential weaknesses of sharp interface models which are related to a loss of robustness. Moreover, a phase field approach provides flexibility in data-driven model order reduction for efficient numerical shape optimization. 

To achieve these goals we combine the porous medium approach of \cite{Hinze_BP03} and a phase field approach including a regularization by the Ginzburg--Landau energy. This results in a diffuse interface problem, which approximates a sharp interface problem for shape optimization in fluids that is penalized by a perimeter term. The related optimization problem then is a control in the coefficient optimal control problem where the phase field represents the control. For the fast numerical solution of those optimal control problems we use POD-MOR techniques which are based upon the findings and methods presented in Section 3.
%
%
%

%
%

 \end{document}